\documentclass[11pt]{amsart}
\usepackage{amsfonts}
\usepackage{graphicx}
\usepackage{amscd}
\usepackage{rotating}
\usepackage{amsmath,amsfonts,amssymb,amsthm,mathrsfs,xypic}
\xyoption{curve}
\usepackage{fancybox}

\newtheorem{thm}{Theorem}[section]

\newtheorem{cor}[thm]{Corollary}

\newtheorem{lem}[thm]{Lemma}

\newtheorem{prop}[thm]{Proposition}

\numberwithin{equation}{section}

\newcommand{\ik}{\operatorname {IKer}}
\newcommand{\pc}{\operatorname {PCoker}}

\newcommand{\A}{\Lambda}
\newcommand{\s}{\hfill\blacksquare}

\newcommand{\Ima}{\operatorname{Im}}
\newcommand{\Cok}{\operatorname{Coker}}
\newcommand{\Ker}{\operatorname{Ker}}
\newcommand{\Hom}{\operatorname{Hom}}

\begin{document}
\title [Auslander-Reiten translations in monomorphism categories]{Auslander-Reiten translations \\ in monomorphism categories}
\author [Bao-Lin Xiong\ \ \ \ \ Pu Zhang\ \ \ \ \ Yue-Hui Zhang]
{Bao-Lin Xiong\ \ \ \ \ Pu Zhang$^*$\ \ \ \ \ Yue-Hui Zhang}
\thanks{$^*$The corresponding author.}
\thanks{Supported by the NSF of China (10725104), and STCSM (09XD1402500).}
\address{Department of Mathematics, \ \
Shanghai Jiao Tong University,  \ \ Shanghai 200240, P. R. China}
\address {xiongbaolin$\symbol{64}$gmail.com, pzhang$\symbol{64}$sjtu.edu.cn, zyh$\symbol{64}$sjtu.edu.cn}
\maketitle
\begin{abstract} We generalize Ringel and Schmidmeier's theory on
the Auslander-Reiten translation of the submodule category $\mathcal
S_2(A)$ to the monomorphism category $\mathcal S_n(A)$. As in the
case of $n=2$, $\mathcal S_n(A)$ has Auslander-Reiten sequences, and
the Auslander-Reiten translation $\tau_{\mathcal{S}}$ of $\mathcal
S_n(A)$ can be explicitly formulated via $\tau$ of $A$-mod.
Furthermore, if $A$ is a selfinjective algebra, we study the
periodicity of $\tau_{\mathcal{S}}$ on the objects of $\mathcal
S_n(A)$, and of the Serre functor $F_{\mathcal S}$ on the objects of
the stable monomorphism category $\underline{\mathcal{S}_n(A)}$. In
particular, $\tau_{\mathcal S}^{2m(n+1)}X\cong X$ for
$X\in\mathcal{S}_n(\A(m, t))$; and $F_{\mathcal S}^{m(n+1)}X\cong X$
for $X\in\underline{\mathcal{S}_n(\A(m, t))}$, where $\A(m, t), \
m\ge1, \ t\ge2,$ are the selfinjective Nakayama algebras.

\vskip5pt

\end{abstract}

\vskip10pt

 \centerline {\bf  Introduction}

\vskip10pt

\vskip10pt

Throughout this paper, $n\ge 2$ is an integer, $A$ an Artin algebra,
and $A$-mod the category of finitely generated left $A$-modules. Let
$\mathcal S_n(A)$ denote the monomorphism category of $A$ (it is
usually called the submodule category if $n=2$).

The study of such a category goes back to G. Birkhoff [B], in which
he initiates to classify the indecomposable objects of $\mathcal
S_2(\Bbb Z/\langle p^t\rangle)$ (see also [RW]). In [Ar], $\mathcal
S_n(R)$ is denoted by $\mathcal C(n, R)$, where $R$ is a commutative
uniserial ring; the complete list of $\mathcal C(n, R)$ of finite
type, and of the representation types of $\mathcal C(n, k[x]/\langle
x^t\rangle)$, are given by D. Simson [S] (see also [SW]). Recently,
after the deep and systematic work of C. M. Ringel and M.
Schmidmeier ([RS1] - [RS3]), the monomorphism category receives more
attention. X. W. Chen [C] shows that $\mathcal{S}_2({\mathcal A})$
of a Frobenius abelian category $\mathcal{A}$ is a Frobenius exact
category. D. Kussin, H. Lenzing, and H. Meltzer [KLM] establish a
surprising link between the stable submodule category with the
singularity theory via weighted projective lines of type $(2, 3,
p)$. In [Z], $\mathcal S_n(\mathcal X)$ is studied for any full
subcategory $\mathcal X$ of $A$-mod, and it is proved that for a
cotilting $A$-module $T$, there is a cotilting $T_n(A)$-module ${\bf
m}(T)$ such that $\mathcal S_n(^\perp T) = \ ^\perp {\rm \bf m}(T)$,
where $T_n(A)=\left(\begin{smallmatrix}
             A & A  & \cdots &  A \\
             0 & A  & \cdots &  A \\
             \vdots & \vdots  & \ddots & \vdots \\
             0 & 0  & \cdots & A
           \end{smallmatrix}\right)_{n\times n}$ is the upper triangular matrix
           algebra of $A$, and $^\perp T$ is the left perpendicular
           category of $T$.  As a consequence, for a Gorenstein algebra $A$,
$\mathcal S_n(^\perp A)$ is exactly the category of
Gorenstein-projective $T_n(A)$-modules.

\vskip10pt

Ringel and Schmidmeier construct minimal monomorphisms, and then
prove that $\mathcal S_2(A)$ is functorially finite in $T_2(A)$-mod.
As a result, $\mathcal S_2(A)$ has Auslander-Reiten sequences.
Surprisingly, the Auslander-Reiten translation $\tau_{\mathcal S}$
of $\mathcal S_2(A)$ can be explicitly formulated as $\tau_{\mathcal
S} X \cong {\bf Mimo} \ \tau \ {\bf Cok}X$ for $X\in
\mathcal{S}_2(A)$ ([RS2], Theorem 5.1), where $\tau$ is the
Auslander-Reiten translation of $A$-mod. Applying this to
selfinjective algebras, among others they get $\tau_{\mathcal S}^6
X\cong X$ for indecomposable nonprojective object
$X\in\mathcal{S}_2(A)$, where $A$ is a commutative uniserial
algebra.

\vskip10pt

A beautiful theory should have a general version. The aim of this
paper is to generalize Ringel and Schmidmeier's work on $\mathcal
S_2(A)$  to $\mathcal S_n(A)$. As in the case of $n=2$, $\mathcal
S_n(A)$ has Auslander-Reiten sequences, and $\tau_{\mathcal{S}}$ of
$\mathcal S_n(A)$ can be formulated in the same form as above: these
can be achieved by using the idea in [RS2]. For selfinjective
algebras, Sections 3 and 4 of this paper contain new considerations.
In order to express the higher power of $\tau_{\mathcal{S}}$, we
need the concept of a rotation of an object in ${\rm
Mor}_n(A\mbox{-}\underline {\rm mod})$, which is defined in [RS2]
for $n=2$. In the general case, such a suitable definition needs to
be chosen from different possibilities, and difficulties need to be
overcome to justify that it is well-defined. Also, the Octahedral
Axiom is needed in computing the higher power of the rotations,
which is the key step in studying the periodicity of
$\tau_{\mathcal{S}}$ and the Serre functor on the objects.

\vskip10pt

We outline this paper. In Section 1 we set up some basic properties
of the categories ${\rm Mor}_n(A)$, $\mathcal{S}_n(A)$ and
$\mathcal{F}_n(A)$, and of the functors ${\bf m}_i$,  ${\bf p}_i$,
${\bf Ker}$, ${\bf Cok}$, ${\bf Mono}$, ${\bf Epi}$; and the
construction of ${\bf Mimo}$. Section 2 is to transfer the
Auslander-Reiten sequences of ${\rm Mor}_n(A)$ to those of
$\mathcal{S}_n(A)$; and to give a formula for $\tau_{\mathcal{S}}$
of $\mathcal S_n(A)$ via $\tau$ of $A$-mod (Theorem 2.4).

In Section 3, $A$ is a selfinjective algebra, and hence the stable
category $A\mbox{-}\underline{\rm mod}$ is a trianglated category
([H]), and $\tau$ is a triangle functor of $A\mbox{-}\underline{\rm
mod}$. Using the rotation and the Octahedral Axiom, we get a formula
for $\underline {\tau_{\mathcal S}^{j(n+1)}X}\in {\rm
Mor}_n(A\mbox{-}\underline {\rm mod})$ for $X\in \mathcal S_n(A)$
and $j\ge 1$ (Theorem 3.4). This can be applied to the study of the
periodicity of $\tau_{\mathcal S}$ on objects. In particular,
$\tau_{\mathcal S}^{2m(n+1)}X\cong X$ for  $X\in \mathcal{S}_n(\A(m,
t))$ (Corollary 3.6), where $\A(m, t), \ m\ge1, \ t\ge2,$ are the
selfinjective Nakayama algebras.

In Section 4, $A$ is a finite-dimensional self-injective algebra
over a field. By \cite{Z}, $\mathcal{S}_n(A)$ is exactly the
category of Gorenstein projective $T_n(A)$-modules, and hence the
stable monomorphism category $\underline{\mathcal{S}_n(A)}$ is a
Hom-finite Krull-Schmidt triangulated category with Auslander-Reiten
triangles. By Theorem I.2.4 of I. Reiten and M. Van den Bergh
\cite{RV}, $\underline{\mathcal{S}_n(A)}$ has a Serre functor
$F_{\mathcal{S}}$. We study the periodicity of $F_{\mathcal S}$ on
the objects of $\underline{\mathcal{S}_n(A)}$ (Theorem 4.3). In
particular, $F_{\mathcal S}^{m(n+1)}X\cong X$ for $X\in
\underline{\mathcal{S}_n(\A(m, t))}$ (Corollary 4.4).

In order to make the main clue clearer, we put the proofs of Lemmas
1.5 and 1.6 in Appendix 1.

Note that $\mathcal{S}_{n, 2}$, $\mathcal{S}_{2, 3}$,
$\mathcal{S}_{2, 4}$, $\mathcal{S}_{2, 5}$, $\mathcal{S}_{3, 3}$ and
$\mathcal{S}_{4, 3}$ are the only representation-finite cases among
all $\mathcal{S}_{n, t}=\mathcal{S}_n(k[x]/\langle x^t\rangle), \
n\ge 2, \ t\ge 2$ (\cite{S}, Theorems 5.2 and 5.5). The
Auslander-Reiten quivers of $\mathcal{S}_{2, t}$ with $t=2, 3, 4, 5$
are given in [RS3]. In Appendix 2, we give the remaining cases. We
also include the AR quivers of $\mathcal{S}_n(\A(2, 2))$ with $n=3$
and $4$.

\section{\bf Basics of morphism categories}

We set up some basic properties of several categories and functors,
which will be used throughout this paper.

\subsection{} \  An object of {\it the morphism category}
${\rm Mor}_n(A)$ is $X_{(\phi_i)}=\left(\begin{smallmatrix}
X_1\\
\vdots\\ X_n
\end{smallmatrix}\right)_{(\phi_i)}$, where $\phi_i: X_{i+1}\rightarrow X_i$
are $A$-maps for $1\le i\le n-1$; and a morphism $f = (f_i):
X_{(\phi_i)}\rightarrow Y_{(\psi_i)}$ is
$\left(\begin{smallmatrix} f_1\\
\vdots\\ f_n
\end{smallmatrix}\right)$,
where $f_i: X_i \rightarrow Y_i$ are $A$-maps for $1\le i\le n$,
such that the following diagram commutes \[\xymatrix@R=0.5cm@C=0.8cm
{ X_n\ar[d]_{f_n}\ar[r]^-{\phi_{n-1}} &
X_{n-1}\ar[d]^{f_{n-1}}\ar[r] & \cdots
\ar[r]& X_2\ar[r]^-{\phi_1}\ar[d]^{f_2} & X_1\ar[d]^{f_1}\\
Y_n\ar[r]^-{\psi_{n-1}} & Y_{n-1}\ar[r] & \cdots \ar[r]&
Y_2\ar[r]^-{\psi_1} & Y_1.}\eqno \ (1.1)\] We call $X_i$ the $i$-th
branch of $X_{(\phi_i)}$, and $\phi_i$ the $i$-th morphism of
$X_{(\phi_i)}$. It is well-known that ${\rm Mor}_n(A)$ is equivalent
to $T_n(A)$-mod (see e.g. [Z], 1.4).  Let $Z_{(\theta_i)}\stackrel
f\rightarrow Y_{(\psi_i)}\stackrel g \rightarrow X_{(\phi_i)}$ be a
sequence in ${\rm Mor}_n(A)$. Then it is exact at $Y_{(\psi_i)}$ if
and only if each sequence $Z_i\stackrel {f_i}\rightarrow
Y_i\stackrel {g_i} \rightarrow X_i$ in $A$-mod is exact at $Y_i$ for
each $1\le i\le n$.

\vskip10pt

By definition, {\it the monomorphism category} $\mathcal{S}_n(A)$ is the full subcategory of 
${\rm Mor}_n(A)$ consisting of the objects $X_{(\phi_i)}$, where
$\phi_i: X_{i+1}\rightarrow X_{i}$ are monomorphisms for $1\leq
i\leq n-1$. Dually, {\it the epimorphism category} $\mathcal F_n(A)$
is the full subcategory of ${\rm Mor}_n(A)$ consisting of the
objects $X_{(\phi_i)}$, where $\phi_i: X_{i+1}\rightarrow X_{i}$ are
epimorphisms for $1\leq i\leq n-1$. Since $\mathcal{S}_n(A)$ and
$\mathcal{F}_n(A)$ are closed under direct summands and extensions,
it follows that they are exact Krull-Schmidt categories, with the
exact structure in ${\rm Mor}_n(A)$.

\vskip10pt

{\it The kernel functor} ${\bf Ker}: {\rm
Mor}_n(A)\rightarrow\mathcal{S}_n(A)$ is given by
$$\left(\begin{smallmatrix} X_1\\ X_2 \\ \vdots\\
X_{n-1} \\ X_n
\end{smallmatrix}\right)_{(\phi_i)} \mapsto \left(\begin{smallmatrix}
X_n\\ \Ker(\phi_1\cdots\phi_{n-1}) \\ \vdots\\ \Ker(\phi_{n-2}\phi_{n-1}) \\
\Ker\phi_{n-1}
\end{smallmatrix}\right)_{(\phi'_i)},$$
where $\phi'_1: \Ker(\phi_1\cdots\phi_{n-1})\hookrightarrow X_n$,
and $\phi'_i: \Ker(\phi_i\cdots\phi_{n-1})\hookrightarrow
\Ker(\phi_{i-1}\cdots\phi_{n-1}), \ 2\le i\le n-1,$ are the
canonical monomorphisms. For a morphism $f: X\rightarrow Y$ in ${\rm
Mor}_n(A)$, ${\bf Ker} f: {\bf Ker}X\rightarrow {\bf Ker}Y$ is
naturally defined via a commutative diagram induced from (1.1). We
also need {\it the cokernel functor} ${\bf Cok}: {\rm
Mor}_n(A)\rightarrow \mathcal{F}_n(A)$ given by
$$\left(\begin{smallmatrix}
X_1\\ X_2 \\ \vdots\\ X_{n-1} \\ X_n
\end{smallmatrix}\right)_{(\phi_i)} \mapsto \left(\begin{smallmatrix}
\Cok\phi_1\\ \Cok(\phi_1\phi_2) \\ \vdots\\ \Cok(\phi_1\cdots\phi_{n-1}) \\
X_1
\end{smallmatrix}\right)_{(\phi''_i)}
,$$ where $\phi''_i: \Cok (\phi_1\cdots
\phi_{i+1})\twoheadrightarrow \Cok (\phi_1\cdots \phi_{i}), \ 1\le
i\le n-2,$ and $\phi''_{n-1}: X_1\twoheadrightarrow
\Cok(\phi_1\cdots \phi_{n-1})$ are the canonical epimorphisms. It is
clear that the restriction of the kernel functor ${\bf Ker}:
\mathcal{F}_n(A)\rightarrow\mathcal{S}_n(A)$ is an equivalence with
quasi-inverse the restriction of the cokernel functor ${\bf Cok}:
\mathcal{S}_n(A)\rightarrow\mathcal{F}_n(A)$.

\subsection{} \ For each $1\leq i\leq n$, the functors ${\bf m}_i:
A$-mod$\rightarrow\mathcal{S}_n(A)$ and ${\bf p}_i:
A$-mod$\rightarrow\mathcal{F}_n(A)$ are defined as follows. For
$M\in A$-mod, $$({\bf m}_i(M))_j = \begin{cases} M, & 1\leq j\leq i; \\
 0, & i+1 \le j\leq n; \end{cases} \ \ \ \ \ ({\bf p}_i(M))_j = \begin{cases} 0, & 1\leq j\leq n-i; \\
 M, & n-i+1\le j\leq n. \end{cases}$$
The $j$-th morphism of ${\bf m}_i(M)$ is ${\rm id}_M$ if $1\leq
j<i$, and $0$ if $i\leq j\leq n-1$; and the $j$-th morphism of ${\bf
p}_i(M)$ is $0$ if $1\leq j < n-i+1$, and ${\rm id}_M$ if $n-i+1\leq
j\leq n-1$. For an $A$-map $f: M\rightarrow N$, define
$${\rm \bf m}_i(f) = \left(\begin{smallmatrix}
f\\
\vdots\\ f \\ 0 \\ \vdots \\ 0
\end{smallmatrix}\right): \ \left(\begin{smallmatrix}
M\\
\vdots\\ M \\ 0 \\ \vdots \\ 0
\end{smallmatrix}\right)\rightarrow  \left(\begin{smallmatrix}
N\\
\vdots\\ N \\ 0 \\ \vdots \\ 0
\end{smallmatrix}\right); \  \ \ \ \ \ {\rm \bf p}_i(f) = \left(\begin{smallmatrix}
0\\
\vdots\\ 0 \\ f \\ \vdots \\ f
\end{smallmatrix}\right): \ \left(\begin{smallmatrix}
0\\
\vdots\\ 0 \\ M \\ \vdots \\ M
\end{smallmatrix}\right)\rightarrow  \left(\begin{smallmatrix}
0\\
\vdots\\ 0 \\ N \\ \vdots \\ N
\end{smallmatrix}\right).$$
 We have  $${\bf Ker}
\ {\bf p}_i(M)={\bf m}_{n-i+1}(M),  \ \ \ \ {\bf Cok} \ {\bf
m}_i(M)={\bf p}_{n-i+1}(M), \ \ \ 1\leq i\leq n.  \eqno \ (1.2)$$

\vskip10pt

\begin{lem}\label{Nakayamafunctor}\quad $(i)$  If $P$ runs over all the indecomposable projective
$A$-modules, then ${\bf m}_1(P), \ \cdots, \ {\bf m}_n(P)$ are the
all indecomposable projective objects in ${\rm Mor}_n(A)$.

\vskip 5pt

$(ii)$ If $I$ runs over all the indecomposable injective
$A$-modules, then ${\bf p}_1(I), \ \cdots, \ {\bf p}_n(I)$ are the
all indecomposable injective objects in ${\rm Mor}_n(A)$.

\vskip 5pt

$(iii)$ The indecomposable projective objects in $\mathcal{S}_n(A)$
are exactly those in ${\rm Mor}_n(A)$.

\vskip 5pt

$(iv)$ If $I$ runs over all the indecomposable injective
$A$-modules, then ${\bf m}_1(I), \ \cdots, \ {\bf m}_n(I)$ are the
all indecomposable injective objects in $\mathcal S_n(A)$.

\vskip5pt

$(v)$ If $P$ runs over all the indecomposable projective
$A$-modules, then ${\bf p}_1(P), \ \cdots, \ {\bf p}_n(P)$ are the
all indecomposable projective objects in $\mathcal F_n(A)$.

\vskip5pt

$(vi)$ The indecomposable injective objects in $\mathcal{F}_n(A)$
are exactly those in ${\rm Mor}_n(A)$.

\vskip5pt

$(vii)$   Let $\mathcal{N}_{\mathcal M}$ and $\mathcal{N}$ be the
Nakayama functor of ${\rm Mor}_n(A)$  and of $A$-mod, respectively.
Then for a projective $A$-module $P$, $\mathcal{N}_{\mathcal M}{\bf
m}_i(P)={\bf p}_{n-i+1}(\mathcal{N}P)$, $1\leq i\leq n$.
\end{lem}
\noindent{\bf Proof.} For convenience, we include a justification.
$(i)$ can be seen from the equivalence ${\rm Mor}_n(A)\cong
T_n(A)$-mod. For $(ii)$, see e.g. Lemma 1.3$(ii)$ in [Z]. $(iii)$
follows from $(i)$, and $(vi)$ follows from $(ii)$. Using the
equivalence ${\bf Ker}: \mathcal{F}_n(A)\rightarrow\mathcal{S}_n(A)$
together with $(vi)$ and (1.2), we see $(iv)$. Using the equivalence
${\bf Cok}: \mathcal{S}_n(A)\rightarrow\mathcal{F}_n(A)$ together
with $(iii)$ and (1.2), we see $(v)$. To see $(vii)$, note that if
$P$ is indecomposable, then $\mathcal{N}_{\mathcal M}{\bf m}_i(P)$
is an indecomposable injective $T_n(A)$-module, hence by $(ii)$ it
is of the form ${\bf p}_j(I)$. Thus
$$\left(\begin{smallmatrix}
0\\
\vdots\\ 0 \\ {\rm soc}(I) \\ 0 \\ \vdots \\ 0
\end{smallmatrix}\right) = {\rm soc} ({\bf p}_j(I)) = {\rm soc} (\mathcal{N}_{\mathcal M}{\bf m}_i(P))
= {\rm top} ({\bf m}_i(P)) = \left(\begin{smallmatrix}
0\\
\vdots\\ 0 \\ {\rm top}(P) \\ 0 \\ \vdots \\ 0
\end{smallmatrix}\right).$$
Thus $n-j+1 = i, \ {\rm soc}(I) = {\rm top}(P)$, which means
$\mathcal{N}_{\mathcal M}{\bf m}_i(P)={\bf
p}_{n-i+1}(\mathcal{N}P)$. $\s$

\vskip10pt

\subsection{}\ Recall the functors ${\bf Mono}: {\rm Mor}_n(A)\rightarrow\mathcal{S}_n(A)$ and
${\bf Epi}: {\rm Mor}_n(A)\rightarrow\mathcal{F}_n(A)$. The first
one is given by
$$\left(\begin{smallmatrix} X_1\\ X_2 \\ \vdots\\
X_{n-1} \\ X_n
\end{smallmatrix}\right)_{(\phi_i)} \mapsto \left(\begin{smallmatrix}
X_1\\ \Ima\phi_1 \\ \vdots\\ \Ima(\phi_1\cdots\phi_{n-2}) \\
\Ima(\phi_1\cdots\phi_{n-1})
\end{smallmatrix}\right)_{(\phi'_i)},$$
where $\phi'_1: \Ima\phi_1\hookrightarrow X_1$, and $\phi'_i:
\Ima(\phi_1\cdots\phi_i)\hookrightarrow
\Ima(\phi_1\cdots\phi_{i-1}), \ 2\le i\le n-1,$ are the canonical
monomorphisms. The second one is given by
$$\left(\begin{smallmatrix}
X_1\\ X_2 \\ \vdots\\ X_{n-1} \\ X_n
\end{smallmatrix}\right)_{(\phi_i)} \mapsto \left(\begin{smallmatrix}
\Ima(\phi_1\cdots\phi_{n-1})\\ \Ima(\phi_2\cdots\phi_{n-1}) \\ \vdots\\ \Ima\phi_{n-1} \\
X_n
\end{smallmatrix}\right)_{(\phi''_i)},$$
where $\phi''_i: \Ima (\phi_{i+1}\cdots
\phi_{n-1})\twoheadrightarrow \Ima (\phi_i\cdots \phi_{n-1}), \ 1\le
i\le n-2,$ and $\phi''_{n-1}: X_n\twoheadrightarrow \Ima\phi_{n-1}$
are the canonical epimorphisms. Then
$${\bf Epi}\cong{\bf Cok} \ {\bf Ker}, \ \ \ {\bf Mono}\cong{\bf Ker} \ {\bf Cok}, \eqno \ (1.3)$$
and hence  ${\bf Mono} X \cong X$ for $X\in \mathcal S_n(A)$, and
${\bf Epi} Y\cong Y$ for $Y\in \mathcal F_n(A)$.

\vskip10pt

For an $A$-map $f: X\rightarrow Y$, denote the canonical $A$-maps
$X\twoheadrightarrow\Ima{f}$ and $\Ima{f}\hookrightarrow Y$ by
$\widetilde{f}$ and ${\rm incl}$, respectively. The following lemma
can be similarly proved as in [RS2] for $n=2$.

\vskip10pt

\begin{lem} {\label{ MonoEpi}}
\ Let $X=X_{(\phi_i)}\in {\rm Mor}_n(A)$. Then \vskip5pt $(i)$ \ The
morphism
$\left(\begin{smallmatrix}1_{X_1}\\\widetilde{\phi_1}\\
\vdots\\\widetilde{\phi_1\cdots\phi_{n-1}}
\end{smallmatrix}\right): \ X\twoheadrightarrow {\bf
Mono} X$ is a left minimal approximation of $X$ in
$\mathcal{S}_n(A)$.

\vskip5pt

$(ii)$ \ The morphism $\left(\begin{smallmatrix}{\rm
incl}\\\vdots\\{\rm incl}\\1_{X_n}
\end{smallmatrix}\right): \ {\bf Epi} X\hookrightarrow X$ is a right minimal
approximation of $X$ in $\mathcal{F}_n(A)$.
\end{lem}

\subsection{}\ For $X_{(\phi_i)}\in{\rm
Mor}_n(A)$, we define ${\bf Mimo}X_{(\phi_i)}\in\mathcal{S}_n(A)$
and ${\bf Mepi}X_{(\phi_i)}\in\mathcal{F}_n(A)$ as follows (see
[Z]). For each $1\leq i\leq n-1,$ fix an injective envelope $e_i':
\Ker\phi_i\hookrightarrow \ik\phi_i$. Then we have an $A$-map $e_i:
X_{i+1}\rightarrow \ik\phi_i$, which is an extension of $e_i'$.
Define ${\bf Mimo}X_{(\phi_i)}$ to be
$$\left(\begin{smallmatrix}
X_1\oplus\ik\phi_{1}\oplus\cdots\oplus\ik\phi_{n-1}\\
X_{2}\oplus\ik\phi_{2}\oplus\cdots\oplus \ik\phi_{n-1}\\ \vdots
\\X_{n-1}\oplus\ik\phi_{n-1}\\X_n
\end{smallmatrix}\right)_{(\theta_i)}, \ \ \ \  {\rm where}\ \ \
 \theta_i = \left(\begin{smallmatrix}
\phi_i& 0 & 0  &\cdots &  0 \\
e_i& 0 & 0 & \cdots &  0 \\
0 & 1& 0 & \cdots &  0 \\
0 & 0& 1 & \cdots &  0\\
 \vdots & \vdots & \vdots &
\ddots & \vdots \\
0 & 0& 0 & \cdots & 1
\end{smallmatrix}\right)_{(n-i+1)\times (n-i)}.$$
By construction  ${\bf Mimo} X_{(\phi_i)}\in \mathcal S_n(A)$. Since
$e_1, \cdots, e_{n-1}$ are not unique, we need to verify that ${\bf
Mimo} X_{(\phi_i)}$ is well-defined. This can be seen from Lemma
1.3$(i)$ below.

 The
object ${\bf Mepi} X_{(\phi_i)}$ is dually defined. Namely, for each
$1\leq i\leq n-1$, fix a projective cover $\pi_i':
\pc\phi_i\twoheadrightarrow\Cok\phi_i$, then we have an $A$-map
$\pi_i: \pc\phi_i\rightarrow X_i$, which is a lift of $\pi_i'$, and
define ${\bf Mepi} X_{(\phi_i)}\in \mathcal{F}_n(A)$ to be
$$\left(\begin{smallmatrix}
X_1\\
X_{2}\oplus\pc\phi_{1}\\\vdots
\\X_{n-1}\oplus\pc\phi_{n-2}\oplus\cdots\oplus\pc\phi_{1}\\X_{n}\oplus\pc\phi_{n-1}\oplus\cdots\oplus\pc\phi_{1}
\end{smallmatrix}\right)_{(\sigma_i)}, \ \ \ \  {\rm where}\ \ \
\sigma_i = \left(\begin{smallmatrix}
\phi_i& \pi_i & 0 & 0&\cdots &  0 \\
0& 0 & 1 & 0 & \cdots &  0 \\
0 & 0& 0 & 1 & \cdots &  0 \\
 \vdots & \vdots & \vdots & \vdots &
\ddots & \vdots\\
0 & 0& 0 & 0 & \cdots & 1
\end{smallmatrix}\right)_{i\times (i+1)}.$$

\vskip10pt

{\bf Remark.} {\it $(i)$ \  ${\bf Mimo} X=X$ for $X\in \mathcal
S_n(A)$,  and ${\bf Mepi} Y = Y$ for $Y\in \mathcal F_n(A)$.

\vskip5pt

$(ii)$ \  If each $X_i$ has no nonzero injective direct summands,
then ${\bf Mimo}X_{(\phi_i)}$ has no nonzero injective direct
summands in $\mathcal{S}_n(A)$. If each $X_i$ has no nonzero
projective direct summands, then ${\bf Mepi}X_{(\phi_i)}$ has no
nonzero projective direct summands in $\mathcal{F}_n(A)$. These can
be seen from Lemma 1.1$(iv)$ and $(v)$, respectively.}

\vskip5pt

\begin{lem}{\label{ MimoMepi}}
\ Let $X\in {\rm Mor}_n(A)$. Then

\vskip5pt

$(i)$ \ The morphism $\left(\begin{smallmatrix}(1, 0, \cdots, 0)\\
\vdots\\(1, 0)\\1\end{smallmatrix}\right):{\bf
Mimo}X\twoheadrightarrow X$ is a right minimal approximation of $X$
in $\mathcal{S}_n(A)$.

\vskip5pt

$(ii)$ \ The morphism
$\left(\begin{smallmatrix}1\\\left(\begin{smallmatrix}1\\0\end{smallmatrix}\right)\\
\vdots\\\left(\begin{smallmatrix}1\\0\\\vdots\\0\end{smallmatrix}\right)\end{smallmatrix}\right):X\hookrightarrow{\bf
Mepi}X$ is a left minimal approximation of $X$ in
$\mathcal{F}_n(A)$.
\end{lem}

For a proof of Lemma 1.3 we refer to [RS2] for $n=2$, and to [Z] in
general case. By Lemmas 1.2 and 1.3, and by Auslander and Smal\o \
[AS], we get the following consequence

\begin{cor} \ {\label{hasAR}}
The subcategories $\mathcal{S}_n(A)$ and $\mathcal{F}_n(A)$ are
functorially finite in ${\rm Mor}_n(A)$ and hence have
Auslander-Reiten sequences.
\end{cor}

This corollary is the starting point of this paper. From now on,
denote by $\tau$, $\tau_{\mathcal M}$, $\tau_{\mathcal S}$ and
$\tau_{\mathcal{F}}$ the Auslander-Reiten translations of $A$-mod,
${\rm Mor}_n(A)$,  $\mathcal{S}_n(A)$ and $\mathcal{F}_n(A)$,
respectively.

\vskip10pt

\subsection{} Let
$A\mbox{-}\underline {\rm mod}$ (resp. $A\mbox{-}\overline {\rm
mod}$) denote the stable category of $A$-mod modulo projective
$A$-modules (resp. injective $A$-modules). Then $\tau = {\rm DTr}$
induces an equivalence $A\mbox{-}\underline {\rm mod}\rightarrow
A\mbox{-}\overline {\rm mod}$ with quasi-inverse $\tau^- = {\rm
TrD}$ ([ARS], p.106). Let ${\rm Mor}_n(A\mbox{-}\overline {\rm
mod})$ denote the morphism category of $A\mbox{-}\overline {\rm
mod}$. Namely, an object of ${\rm Mor}_n(A\mbox{-}\overline {\rm
mod})$ is
$\overline { X_{(\phi_i)}} = X_{( \overline{\phi_i})} = \left(\begin{smallmatrix} X_1\\ \vdots\\
X_n\end{smallmatrix}\right)_{(\overline{\phi_i})}$ with
$\overline{\phi_i}: X_{i+1}\rightarrow X_i$ in $A$-$\overline{\rm
mod}$ for $1\leq i\leq n-1$; and a morphism from
$X_{(\overline{\phi_i})}$
to $Y_{(\overline{\theta_i})}$ is $\left(\begin{smallmatrix} \overline{f_1}\\
\vdots\\
\overline{f_n}\end{smallmatrix}\right)$, such that the corresponding
version of (1.1) commutes in $A$-$\overline{\rm mod}$. Similarly,
one has the morphism category ${\rm Mor}_n(A\mbox{-}\underline {\rm
mod})$, in which an object is denoted by $\underline{ X_{(\phi_i)}}
= X_{(\underline{\phi_i})}$.

\vskip10pt

The following two lemmas will be heavily used in Sections 2 and 3.
In order to make the main clue clearer, we put their proofs in
Appendix 1.

\begin{lem}\label{Mimosummand}
\ Let $X_{(\phi_i)}\in{\rm Mor}_n(A)$.

\vskip5pt

$(i)$ \ Let $I_2, \cdots, I_n$ be injective $A$-modules such
that $X'_{(\phi'_i)}=\left(\begin{smallmatrix}X_1\oplus I_2\oplus\cdots\oplus I_n\\
\vdots\\X_{n-1}\oplus
I_n\\X_n\end{smallmatrix}\right)_{(\phi'_i)}\in \mathcal{S}_n(A)$,
where each $\phi'_i$ is of the form $\left(\begin{smallmatrix}\phi_i
&
* & \cdots & *\\ * & * & \cdots & * \\ \vdots & \vdots & \ddots & \vdots\\ * & * & \cdots & *
\end{smallmatrix}\right)_{(n-i+1)\times(n-i)}.$
Then $X'_{(\phi'_i)}\cong{\bf Mimo}X_{(\phi_i)}\oplus J$, where $J$
is an injective object of $\mathcal{S}_n(A)$.

\vskip5pt

$(ii)$ \ Let $P_1, \cdots, P_{n-1}$ be projective $A$-modules
such that $X''_{(\phi''_i)}=\left(\begin{smallmatrix}X_1\\
X_{2}\oplus P_1 \\ \vdots\\X_n\oplus P_{n-1}\oplus\cdots\oplus
P_1\end{smallmatrix}\right)_{(\phi''_i)}\in \mathcal{F}_n(A)$, where
each $\phi''_i$ is of the form $\left(\begin{smallmatrix}\phi_i &
* & \cdots & *\\ * & * & \cdots & * \\ \vdots & \vdots & \ddots & \vdots\\ * & * & \cdots & *
\end{smallmatrix}\right)_{i\times(i+1)}.$
Then $X''_{(\phi''_i)}\cong{\bf Mepi}X_{(\phi_i)}\oplus L$, where
$L$ is a projective object of $\mathcal{F}_n(A)$.
\end{lem}

\begin{lem}\label{MimoIso} \ Let $X_{(\phi_i)}, \ Y_{(\psi_i)}\in {\rm Mor}_n(A)$.

\vskip5pt

$(i)$  \ \  If  all branches $X_i$ and $Y_i$ have no nonzero
injective direct summands, then ${\bf Mimo}X_{(\phi_i)} \cong {\bf
Mimo} Y_{(\psi_i)}$ in $\mathcal{S}_n(A)$ if and only if
$X_{(\overline{\phi_i})}\cong Y_{(\overline{\psi_i})}$ in ${\rm
Mor}_n(A\mbox{-}\overline {\rm mod})$.

\vskip5pt

$(ii)$ \ If  all $X_i$ and $Y_i$ have no nonzero projective direct
summands, then ${\bf Mepi}X_{(\phi_i)}\cong {\bf Mepi} Y_{(\psi_i)}$
in $\mathcal{F}_n(A)$ if and only if $X_{(\underline{\phi_i})}\cong
Y_{(\underline {\psi_i})}$ in ${\rm Mor}_n(A\mbox{-}\underline {\rm
mod})$. \end{lem}

\vskip10pt

\section{\bf The Auslander-Reiten translation of
$\mathcal{S}_n(A)$}

In this section, we first transfer the Auslander-Reiten sequences of
${\rm Mor}_n(A)$ to those of $\mathcal{S}_n(A)$ and
$\mathcal{F}_n(A)$; and then give a formula of the Auslander-Reiten
translation $\tau_{\mathcal{S}}$ of $\mathcal{S}_n(A)$ via $\tau$ of
$A$-mod. Results and methods in this section are generalizations of
the corresponding ones in the case of $n=2$, due to Ringel and
Schmidmeier [RS2].

\subsection{}  The following fact is crucial for later use.

\begin{lem}{\label{KerCok}}
\ Let \ $0\rightarrow X_{(\phi_i)}\overset{f}\rightarrow
Y\overset{g}\rightarrow Z\rightarrow 0$ be an Auslander-Reiten
sequence of ${\rm Mor}_n(A)$.

\vskip5pt

$(i)$ If \ ${\bf Ker}Z$ is not projective, then $0\rightarrow {\bf
Ker}X\overset{{\bf Ker}f}\rightarrow {\bf Ker}Y\overset{{\bf
Ker}g}\rightarrow {\bf Ker}Z\rightarrow0$ is either split exact, or
an Auslander-Reiten sequence of $\mathcal{S}_n(A)$.

\vskip5pt

$(ii)$ If \ ${\bf Cok}X$ is not injective, then $0\rightarrow {\bf
Cok}X\overset{{\bf Cok}f}\rightarrow {\bf Cok}Y\overset{{\bf
Cok}g}\rightarrow {\bf Cok}Z\rightarrow0$ is either split exact, or
an Auslander-Reiten sequence of $\mathcal{F}_n(A)$.
\end{lem}
\noindent{\bf Proof.} We only prove $(i)$. Put $g'= {\bf Ker} g$ and
$f'={\bf Ker}f$. By Snake Lemma, $0\rightarrow {\bf
Ker}X\overset{f'}\rightarrow {\bf Ker}Y\overset{g'}\rightarrow {\bf
Ker}Z$ is exact. Assume that $g'$ is not a split epimorphism. We
claim that $g'$ is right almost split. Let $v: W\rightarrow {\bf
Ker}Z$ be a morphism in $\mathcal S_n(A)$ which is not a split
epimorphism. Applying ${\bf Cok}$, we get $t'={\bf Cok}v: {\bf Cok}W
\rightarrow {\bf Cok} \ {\bf Ker}Z = {\bf Epi}Z$, which is not a
split epimorphism, and hence the composition $t: {\bf Cok}W\stackrel
{t'}\rightarrow {\bf Epi}Z\stackrel\sigma\hookrightarrow Z$ is not a
split epimorphism. So, there is a morphism $s: {\bf Cok}W\rightarrow
Y$ such that $t=gs$. Applying ${\bf Ker}$, we get ${\bf Ker}t =
g'{\bf Ker}s$ with ${\bf Ker}s: W\rightarrow{\bf Ker}Y$. Since ${\bf
Ker}\sigma = {\rm id}_{{\bf Ker}Z}$, we see $v = g'{\bf Ker}s$. This
proves the claim.

Since $g'$ is right almost split and ${\bf Ker}Z$ is not projective,
it follows that $g'$ is epic, and hence $f'$ is not a split
monomorphism. We claim that $f'$ is left almost split. For this, let
$p: {\bf Ker}X\rightarrow B$ be a morphism in $\mathcal{S}_n(A)$
which is not a split monomorphism. Take an injective envelope
$(e_i): \left(\begin{smallmatrix}
B_1/B_2\\ \vdots\\
B_1/B_n \end{smallmatrix}\right)_{(\pi_i)} \hookrightarrow
\left(\begin{smallmatrix}
I_1\\ \vdots\\
I_{n-1}
\end{smallmatrix}\right)_{(b'_i)}$  in ${\rm
Mor}_{n-1}(A)$. Put $B' = \left(\begin{smallmatrix}
I_1\\ \vdots\\
I_{n-1} \\ B_1
\end{smallmatrix}\right)_{(b'_i)} \in {\rm Mor}_{n}(A)$, where $b'_{n-1}$ is
the composition $B_1\stackrel {\pi_{n-1}}\twoheadrightarrow
B_1/B_n\stackrel {e_{n-1}}\hookrightarrow I_{n-1}$. By construction
we get a morphism $\left(\begin{smallmatrix}
e_1\\ \vdots\\
e_{n-1} \\ {\rm id}_{B_1}
\end{smallmatrix}\right): {\bf Cok}B \rightarrow B'$, and ${\bf Cok}B = {\bf Epi}B'$. Hence
${\bf Ker} B' = B$. Put $r'=(r'_i) = {\bf Cok}p: {\bf
Epi}X\rightarrow {\bf Cok}B$. Then we have a morphism
$\left(\begin{smallmatrix}
e_1r_1'\\ \vdots\\
e_{n-1}r'_{n-1}
\end{smallmatrix}\right): \left(\begin{smallmatrix}
\Ima(\phi_1\cdots\phi_{n-1})\\\Ima(\phi_2\cdots\phi_{n-1})\\ \vdots\\
\Ima\phi_{n-1}
\end{smallmatrix}\right)\rightarrow\left(\begin{smallmatrix}
I_1\\ \vdots\\
I_{n-1}
\end{smallmatrix}\right)_{(b'_i)}$ in ${\rm
Mor}_{n-1}(A)$. By the canonical monomorphism $(j_i):
\left(\begin{smallmatrix}
\Ima(\phi_1\cdots\phi_{n-1})\\\Ima(\phi_2\cdots\phi_{n-1})\\ \vdots\\
\Ima\phi_{n-1}
\end{smallmatrix}\right)\hookrightarrow \left(\begin{smallmatrix}
X_1\\X_2\\ \vdots\\
X_{n-1}
\end{smallmatrix}\right)$ in ${\rm Mor}_{n-1}(A)$, we get a morphism $\left(\begin{smallmatrix}
r_1\\ \vdots\\
r_{n-1}
\end{smallmatrix}\right): \ \left(\begin{smallmatrix}
X_1\\ \vdots\\
X_{n-1}
\end{smallmatrix}\right)_{(\phi_i)}\rightarrow \left(\begin{smallmatrix}
I_1\\ \vdots\\
I_{n-1}
\end{smallmatrix}\right)_{(b'_i)}$ in ${\rm Mor}_{n-1}(A)$, such that $e_ir'_i
= r_ij_i$ for $1\le i\le n-1$. Then we get a morphism $r=(r_i):
X\rightarrow B'$ in ${\rm Mor}_{n}(A)$ by letting $r_n = r'_n$ (we
only need to check  $b'_{n-1}r_n = r_{n-1}\phi_{n-1}$: $b'_{n-1}r_n
= b'_{n-1}r'_n = e_{n-1}\pi_{n-1}r'_n = e_{n-1}r'_{n-1}\phi'_{n-1} =
r_{n-1}j_{n-1}\phi'_{n-1} = r_{n-1}\phi_{n-1}$). By construction we
get  ${\bf Epi} r = r': {\bf Epi}X \rightarrow {\bf Epi}B' = {\bf
Cok}B$. This process can be figured as follows.
$$\xymatrix@R=0.3cm@C=0.5cm{ & X_n\ar@{->>}[dd]^-<>>>{\phi'_{n-1}}\ar[dl]_{{\rm id}_{X_n}}\ar[r]^{r'_n}& B_{1}\ar@{->>}[dd]^<>>>>>{\pi_{n-1}}\ar[rr]^-{{\rm id}_{B_1}}& &B_{1}\ar[dd]^-{b'_{n-1}} \\
X_n\ar[urrrr]_-{r_n =r'_n}\ar[dd]_{\phi_{n-1}} & & &\\
 & \Ima\phi_{n-1}\ar[r]^-{r'_{n-1}}\ar@{->>}[dd]^-<>>>{\phi'_{n-2}}\ar@{_(->}[dl]_-{j_{n-1}} &
 B_1/B_n\ar@{->>}[dd]^-{\pi_{n-2}}\ar@{^(->}[rr]^-{e_{n-1}} & &I_{n-1}\ar[dd]^-{b'_{n-2}} \\
 X_{n-1}\ar@{.>}[urrrr]_{r_{n-1}}\ar[d]
& & & &\\
 \vdots\ar[dd] & \vdots\ar@{->>}[d] & \vdots\ar@{->>}[d] & &\vdots\ar[d]^-{b'_{2}} \\
 & \Ima(\phi_{2}\cdots\phi_{n-1})\ar[r]^-{r'_{2}}\ar@{->>}[dd]^-<>>>{\phi'_{1}}\ar@{_(->}[dl]_-{j_{2}}
& B_1/B_3\ar@{->>}[dd]^{\pi_1}\ar@{^(->}[rr]^-{e_{2}} & &
I_{2}\ar[dd]^-{b'_{1}} \\
 X_{2}\ar[dd]_{\phi_{1}}\ar@{.>}[urrrr]_{r_{2}} & & & &\\
 & \Ima(\phi_{1}\cdots\phi_{n-1})\ar[r]^-{r'_{1}}\ar@{_(->}[dl]_-{j_{1}} & B_1/B_{2}\ar@{^(->}[rr]^-{e_{1}} & &I_{1} \\
 X_{1}\ar@{.>}[urrrr]_{r_{1}} & & &&}$$
Clearly, $r$ is not a spit monomorphism (otherwise, ${\bf Epi} r =
r': {\bf Epi}X \rightarrow {\bf Epi}B' = {\bf Cok}B$ is a spit
monomorphism, and hence $p = {\bf Ker} \ r'$ is a spit
monomorphism). So there is a morphism $h: Y\rightarrow B'$ such that
$r=hf: X\rightarrow B'$. Applying ${\bf Ker}$, we get ${\bf Ker} r =
({\bf Ker} h)  f',$ where ${\bf Ker} h: {\bf Ker} Y\rightarrow {\bf
Ker} B' = B.$  Since ${\bf Ker} r = {\bf Ker} \ {\bf Epi} r = {\bf
Ker} r' = {\bf Ker} \ {\bf Cok} p = p,$ we get $p = ({\bf Ker} h)
f'.$ This proves the claim, and completes the proof. $\s$

\vskip10pt

\begin{prop}\label{NewAR} \ Let $0\rightarrow  X \stackrel f \rightarrow
Y \stackrel g\rightarrow Z\rightarrow0$ be an Auslander-Reiten
sequence of ${\rm Mor}_n(A)$.

\vskip5pt

$(i)$ \ If $Z\in \mathcal{F}_n(A)$, and $Z$ is not projective in
$\mathcal{F}_n(A)$, then
$$
0\rightarrow{\bf Epi}X\stackrel {{\bf Epi}f}\rightarrow {\bf Epi}Y
\stackrel {{\bf Epi}g} \rightarrow Z\rightarrow 0
\label{Episeq}\eqno (2.1)$$ is an Auslander-Reiten sequence of
$\mathcal{F}_n(A)$; and
$${\label{Kerseq}} 0\rightarrow {\bf Ker}X \stackrel {{\bf Ker}f} \rightarrow{\bf
Ker}Y\stackrel {{\bf Ker}g}\rightarrow {\bf Ker}Z\rightarrow0 \eqno
(2.2)$$
 is an Auslander-Reiten sequence of $\mathcal{S}_n(A)$.

\vskip5pt

$(ii)$ \ If $X\in \mathcal{S}_n(A)$, and $X$ is not injective in
$\mathcal{S}_n(A)$, then
$$0\rightarrow X \stackrel {{\bf Mono}f} \longrightarrow{\bf Mono}Y\stackrel {{\bf Mono}g} \longrightarrow {\bf
Mono}Z\rightarrow0 \label{Monoseq}$$
 is an Auslander-Reiten sequence of $\mathcal{S}_n(A)$; and
$${\label{Cokseq}} 0\rightarrow {\bf Cok}X\stackrel {{\bf Cok}f}\rightarrow{\bf
Cok}Y\stackrel {{\bf Cok}g}\rightarrow {\bf Cok}Z\rightarrow 0
$$ is an Auslander-Reiten sequence of $\mathcal{F}_n(A)$.
\end{prop}
\noindent{\bf Proof.} We only show $(i)$. Since ${\bf Ker}:
\mathcal{F}_n(A)\rightarrow\mathcal{S}_n(A)$ is an equivalence, and
$Z$ is not projective in $\mathcal{F}_n(A)$, it follows that ${\bf
Ker}Z$ is not projective, and hence by Lemma \ref{KerCok}, (2.2) is
either split exact, or an Auslander-Reiten sequence in
$\mathcal{S}_n(A)$. Applying the equivalence ${\bf Cok}:
\mathcal{S}_n(A)\rightarrow \mathcal{F}_n(A)$, and using ${\bf
Epi}\cong{\bf Cok} \ {\bf Ker}$ and $Z\in\mathcal{F}_n(A)$, we see
that (2.1) is either split exact, or an Auslander-Reiten sequence in
$\mathcal{F}_n(A)$. While ${\bf Epi}Y\rightarrow Z$ is the
composition of the canonical monomorphism ${\bf Epi}Y\hookrightarrow
Y$ and the right almost split morphism $Y\rightarrow Z$, so ${\bf
Epi}Y\rightarrow Z$ is not a split epimorphism, hence (2.1) is an
Auslander-Reiten sequence in $\mathcal{F}_n(A)$, so is $(2.2)$. $\s$

\vskip 10pt

\subsection{} We have the following
relationship between $\tau_{\mathcal{S}}$ and $\tau_{\mathcal M}$.

\vskip10pt

\begin{cor}\label{TauRelation} \ $(i)$ \ If $Z\in \mathcal{S}_n(A)$, then
$\tau_{\mathcal{S}}Z \cong {\bf Ker} \ \tau_{\mathcal M} \ {\bf
Cok}Z$, and $\tau_{\mathcal{S}}^-Z\cong {\bf Mono}\ \tau_{\mathcal
M}^-Z$.

\vskip 10pt

$(ii)$ \ If $Z\in \mathcal{F}_n(A)$, then $\tau_{\mathcal{F}}Z \cong
{\bf Epi}\ \tau_{\mathcal M} Z$, and $\tau_{\mathcal{F}}^-Z \cong
{\bf Cok}\ \tau_{\mathcal M}^- \ {\bf Ker}Z$.
\end{cor}

\noindent{\bf Proof.}  We only prove the first formula of $(i)$.
Assume that $Z$ is indecomposable. If $Z$ is projective, then  $Z =
{\bf m}_i(P)$ by Lemma 1.1$(iii)$, where $P$ is an indecomposable
projective $A$-module. By the definition of $\tau_{\mathcal M}$ and
a direct computation, we have $\tau_{\mathcal M} {\bf Cok} Z =
\tau_{\mathcal M}{\bf Cok} \ {\bf m}_i(P) = \tau_{\mathcal M} \ {\bf
p}_{n-i+1}(P)
= \left(\begin{smallmatrix} *\\ \vdots\\ * \\
0\end{smallmatrix}\right) $, it follows that ${\bf Ker} \
\tau_{\mathcal M} \ {\bf Cok} Z = 0 = \tau_{\mathcal{S}}Z$. Assume
that $Z\in \mathcal{S}_n(A)$ is not projective. Since ${\bf Cok}:
\mathcal{S}_n(A)\rightarrow \mathcal{F}_n(A)$ is an equivalence,
${\bf Cok}Z\in \mathcal{F}_n(A)$ is not projective. By Lemma
1.1$(i)$ and $(v)$, ${\bf Cok}Z$ is an indecomposable nonprojective
object in ${\rm Mor}_n(A)$. Replacing $Z$ by ${\bf Cok}Z$ in (2.2),
we get the assertion by ${\bf Ker} \ {\bf Cok}Z\cong Z$. $\s$

\vskip 10pt

\subsection{Example.} \ Let $k$ be a field,
$A=k[x]/\langle x^2\rangle$, and $S$ be the simple $A$-module.
Denote by $i: S\hookrightarrow A$ and $\pi: A\twoheadrightarrow S$
the canonical $A$-maps. Then we have the Auslander-Reiten sequence
in ${\rm Mor}_3(A)$
$$\xymatrix@1{0\ar[r] & {\left(\begin{smallmatrix}0\\ A\\S\end{smallmatrix}\right)}_{(i,
0)}\ar[rr]^-{\left(\begin{smallmatrix}\left(\begin{smallmatrix}
0\\0\\1
\end{smallmatrix}\right)\\\left(\begin{smallmatrix}0\\\left(\begin{smallmatrix}\pi\\1
\end{smallmatrix}\right)\\i
\end{smallmatrix}\right)
\end{smallmatrix}\right)}
& & {\left(\begin{smallmatrix}0\\0\\S\end{smallmatrix}\right)}\oplus{\left(\begin{smallmatrix}S\\
S\oplus
A\\A\end{smallmatrix}\right)}_{\left(\left(\begin{smallmatrix}0\\1\end{smallmatrix}\right),(-1,
\pi)\right)}\ar[rr]^-{\left(\left(\begin{smallmatrix}0\\0\\-i
\end{smallmatrix}\right), \left(\begin{smallmatrix}1\\(-1, \pi)\\1
\end{smallmatrix}\right)\right)}& &
\ \
{\left(\begin{smallmatrix}S\\S\\A\end{smallmatrix}\right)}_{(\pi,
1)}\ar[r]& 0 }.$$ By $(2.1)$, we get an Auslander-Reiten sequence in
$\mathcal{F}_3(A)$
$$\xymatrix@1{0\ar[r] & {\left(\begin{smallmatrix}0\\ S\\S\end{smallmatrix}\right)}_{(1, 0)}
\ar[rr]^-{\left(\begin{smallmatrix}\left(\begin{smallmatrix} 0\\0\\1
\end{smallmatrix}\right)\\\left(\begin{smallmatrix} 0 \\ i\\ i
\end{smallmatrix}\right)
\end{smallmatrix}\right)}
& & {\left(\begin{smallmatrix}0\\0\\S\end{smallmatrix}\right)}\oplus{\left(\begin{smallmatrix}S\\
A\\A\end{smallmatrix}\right)}_{(1,
\pi)}\ar[rr]^-{\left(\left(\begin{smallmatrix}0\\0\\-i
\end{smallmatrix}\right), \left(\begin{smallmatrix}1\\\pi\\1
\end{smallmatrix}\right)\right)}& &
{\left(\begin{smallmatrix}S\\S\\A\end{smallmatrix}\right)}_{(\pi,
1)}\ar[r]& 0 }.$$By $(2.2)$, we get an Auslander-Reiten sequence in
$\mathcal{S}_3(A)$
$$\xymatrix@1{0\ar[r] & {\left(\begin{smallmatrix}S\\ S\\0\end{smallmatrix}\right)}_{(0, 1)}
\ar[rr]^-{\left(\begin{smallmatrix}\left(\begin{smallmatrix} 1\\1\\0
\end{smallmatrix}\right)\\\left(\begin{smallmatrix} i \\ 1\\ 0
\end{smallmatrix}\right)
\end{smallmatrix}\right)}
& & {\left(\begin{smallmatrix}S\\S\\S\end{smallmatrix}\right)}_{(1, 1)}\oplus{\left(\begin{smallmatrix}A\\
S\\0\end{smallmatrix}\right)}_{(0,
i)}\ar[rr]^-{\left(\left(\begin{smallmatrix}i\\1\\1
\end{smallmatrix}\right), \left(\begin{smallmatrix}-1\\-1\\0
\end{smallmatrix}\right)\right)}& &
{\left(\begin{smallmatrix}A\\S\\S\end{smallmatrix}\right)}_{(1,
i)}\ar[r]& 0 }.$$

\subsection{} In Corollary \ref {TauRelation}$(i)$, $\tau_{\mathcal{S}}$ is
formulated via $\tau_{\mathcal M}$ of ${\rm Mor}_n(A)$. However,
$\tau_{\mathcal M}$ is usually more complicated than $\tau$. The
rest of this section is to give a formula of $\tau_{\mathcal{S}}$
via $\tau$.

\vskip10pt

Before stating the main result, we need a notation. For
$X_{(\underline{\phi_i})}\in {\rm Mor}_n(A\mbox{-}\underline {\rm
mod})$, define $\tau X_{(\underline {\phi_i})}
=\left(\begin{smallmatrix} \tau X_1\\ \vdots\\
\tau X_n\end{smallmatrix}\right)_{(\tau \underline {\phi_i})}\in
{\rm Mor}_n(A\mbox{-}\overline {\rm mod}).$ Consider the full
subcategory given by
$$\{Y_{(\psi_i)} = \left(\begin{smallmatrix} \tau X_1\\ \vdots\\
\tau X_n\end{smallmatrix}\right)_{(\psi_i)}\in {\rm Mor}_n(A) \ | \
Y_{(\overline {\psi_i})} \cong  \tau X_{(\underline {\phi_i})}\}.$$
Any object in this full subcategory will be denoted by $\tau
X_{({\phi_i})}$ (we emphasize that this convention will cause no
confusions). So, for $X_{(\phi_i)}\in {\rm Mor}_n(A)$ we have
$\overline {\tau X_{({\phi_i})}} \cong \tau X_{(\underline
{\phi_i})}.$ By Lemma 1.6$(i)$, ${\bf Mimo} \ \tau X_{(\phi_i)}$ is
a well-defined object in $\mathcal{S}_n(A)$, and there are
isomorphisms $\overline {{\bf Mimo} \ \tau X_{(\phi_i)}}\cong
\overline {\tau X_{({\phi_i})}} \cong \tau X_{(\underline
{\phi_i})}$ in ${\rm Mor}_n(A\mbox{-}\overline {\rm mod}).$ If $A$
is selfinjective, then ${\rm Mor}_n(A\mbox{-}\overline {\rm mod})=
{\rm Mor}_n(A\mbox{-}\underline {\rm mod}),$ so the isomorphism
above is read as follows, which is needed in the next section
$$\underline {{\bf Mimo} \ \tau X_{(\phi_i)}}\cong \underline {\tau X_{({\phi_i})}} \cong \tau X_{(\underline
{\phi_i})}.\eqno(2.3)$$

\vskip10pt

Similarly, one has the convention $\tau^- \ X_{(\phi_i)}$, and ${\bf
Mepi} \ \tau^- \ X_{(\phi_i)}\in \mathcal{F}_n(A)$ is well-defined.

\vskip10pt The following result is a generalization of Theorem 5.1
of Ringel and Schmidmeier [RS2].

\vskip10pt

\begin{thm}\label{tauformula} Let $X_{(\phi_i)}\in
\mathcal{S}_n(A)$. Then

\vskip5pt

$(i)$ \ \ $\tau_{\mathcal{S}}X_{(\phi_i)}\cong {\bf Mimo} \ \tau \
{\bf Cok}X_{(\phi_i)}.$

\vskip5pt

$(ii)$ \ \ $\tau^{-}_{\mathcal{S}}X_{(\phi_i)} \cong {\bf Ker} \
{\bf Mepi} \ \tau^- \ X_{(\phi_i)}.$
\end{thm}
\noindent {\bf Proof.} We only prove $(i)$. Recall ${\bf
Cok}X_{(\phi_i)} = \left(\begin{smallmatrix}
\Cok\phi_1\\ \Cok(\phi_1\phi_2) \\ \vdots\\ \Cok(\phi_1\cdots\phi_{n-1}) \\
X_1
\end{smallmatrix}\right)_{(\phi'_i)}.$ Fix a minimal projective presentation
$Q_n\stackrel {d_n} \rightarrow P\stackrel {e}\rightarrow
X_1\rightarrow 0$. Then we get the following commutative diagram
with exact rows
\[\xymatrix@R=0.4cm@C=0.4cm {Q_n\ar[rr]^-{d_n}\ar@{.>}[d]_{s_{n-1}} & & P\ar[rr]^-e\ar@{=}[d] & & X_1\ar[r]\ar@{->>}[d]^{\phi'_{n-1}} & 0  \\
Q_{n-1}\ar[rr]^-{d_{n-1}}\ar@{.>}[d] & &
P\ar[rr]^-{\phi'_{n-1}e}\ar@{=}[d] & &
\Cok(\phi_1\cdots\phi_{n-1})\ar[r]\ar@{->>}[d] & 0 \\
\vdots\ar@{.>}[d]_{s_1} & & \vdots\ar@{=}[d] & & \vdots\ar@{->>}[d]^{\phi'_1}\\
Q_1\ar[rr]^-{d_1} & & P\ar[rr]^-{\phi'_1\cdots\phi'_{n-1}e} & &
\Cok(\phi_1)\ar[r]& 0 } \eqno(2.4)\] where $Q_i\twoheadrightarrow
\Ker(\phi'_i\cdots\phi'_{n-1}e)$ is a projective cover, and $d_i$ is
the composition $Q_i\twoheadrightarrow
\Ker(\phi'_i\cdots\phi'_{n-1}e)\hookrightarrow P$. Applying the
Nakayama functor $\mathcal{N}=D\Hom_A(-, \ _AA)$,  we get the
following commutative diagram
$$\xymatrix@R=0.4cm@C=0.4cm{0\ar[r] & \tau X_1\ar[rr]^-{\sigma_n}\ar@{.>}[d]_-{\alpha_{n-1}} & &
\mathcal{N}Q_n\ar[rr]^-{\mathcal{N}d_n}\ar[d]_{\mathcal{N}s_{n-1}} &
&
\mathcal{N}P\ar@{=}[d]\\
0\ar[r] &
\tau\Cok(\phi_1\cdots\phi_{n-1})\ar[rr]^-{\sigma_{n-1}}\ar@{.>}[d]_-{\alpha_{n-2}}
& &
\mathcal{N}Q_{n-1}\ar[rr]^-{\mathcal{N}d_{n-1}}\ar[d]_{\mathcal{N}s_{n-2}}
& &
\mathcal{N}P\ar@{=}[d]\\
 & \vdots\ar@{.>}[d] & & \vdots\ar[d] & & \vdots\ar@{=}[d]\\
  0\ar[r] & \tau \ \Cok(\phi_1\phi_2)\ar[rr]^-{\sigma_2}\ar@{.>}[d]_-{\alpha_1}& &
\mathcal{N}Q_2\ar[rr]^-{\mathcal{N}d_2}\ar[d]_-{\mathcal{N}s_1} & & \mathcal{N}P\ar@{=}[d]\\
 0\ar[r] & \tau \ \Cok\phi_1\ar[rr]^-{\sigma_1}& &
\mathcal{N}Q_1\ar[rr]^-{\mathcal{N}d_1} & & \mathcal{N}P.
}\eqno(2.5)$$

{\bf Step 1.} By $(2.4)$, we get a projective presentation
$$\xymatrix@1{\bigoplus\limits_{i=1}^n{\bf
m}_i(Q_i)\ar[rrrr]^-{\left(\left(\begin{smallmatrix}
d_1 \\0\\0\\
\vdots\\0
\end{smallmatrix}\right), \left(\begin{smallmatrix}
d_2 \\d_2\\0\\
\vdots\\0
\end{smallmatrix}\right), \cdots, \left(\begin{smallmatrix}
d_n \\d_n\\d_n\\
\vdots\\d_n
\end{smallmatrix}\right)\right)} & & & &{\bf
m}_n(P)\ar[rr]^-{\left(\begin{smallmatrix}
\phi'_1\cdots\phi'_{n-1}e\\
\vdots\\
\phi'_{n-1}e\\
e
\end{smallmatrix}\right)} & &
{\bf Cok}X_{(\phi_i)}\ar[r] & 0}\eqno(2.6)$$ (the exactness can be
seen as follows: by $(2.4)$ we have $\Ima d_n\subseteq \Ima
d_{n-1}\subseteq\cdots\subseteq\Ima d_1$, and hence
$Q_i\oplus\cdots\oplus Q_n\stackrel{(d_i, \cdots,
d_n)}\longrightarrow
P\stackrel{\phi'_i\cdots\phi'_{n-1}e}\longrightarrow\Cok(\phi_i\cdots\phi_{n-1})\rightarrow0$
is exact). In order to obtain a minimal projective presentation from
$(2.6)$, we have to split off a direct summand of
$\bigoplus\limits_{i=1}^n{\bf m}_i(Q_i)$. By Lemma 1.1$(i)$, this
direct summand is of the form $\bigoplus\limits_{i=1}^{n-1} {\bf
m}_i(Q'_i)$ where $Q'_i$ is a direct summand of $Q_i, \ 1\leq i\leq
n-1$,  since $\left(\begin{smallmatrix}
\phi'_1\cdots\phi'_{n-1}e\\
\vdots\\
\phi'_{n-1}e\\
e
\end{smallmatrix}\right)$ is already
minimal and $Q_n\stackrel {d_n} \rightarrow P\stackrel
{e}\rightarrow X_1\rightarrow 0$ is already a minimal projective
presentation. Applying the Nakayama functor $\mathcal{N}_{\mathcal
M}$, we get the exact sequence
$$\xymatrix@1{0\ar[r] & \tau_{\mathcal M} {\bf Cok}X_{(\phi_i)}\oplus
\mathcal{N}_{\mathcal M}(\bigoplus\limits_{i=1}^{n-1}{\bf m}_i(Q'_i))\ar[r]&
\mathcal{N}_{\mathcal M}(\bigoplus\limits_{i=1}^n {\bf
m}_i(Q_i))\ar[r]& \mathcal{N}_{\mathcal M}{\bf m}_n(P).}$$ By Lemma
1.1$(vii)$ this exact sequence can be written as
$$\xymatrix@1{0\ar[r] & \tau_{\mathcal M} {\bf
Cok}X_{(\phi_i)}\oplus\bigoplus\limits_{i=1}^{n-1}{\bf
p}_{n-i+1}(\mathcal{N}Q'_{i})\ar[r]& \bigoplus\limits_{i=1}^n{\bf
p}_{n-i+1}(\mathcal{N}Q_{i})\ar[r]^-d& {\bf p}_1(\mathcal{N}P)}
\eqno (2.7)$$where
$d=\left(\left(\begin{smallmatrix} 0 \\
\vdots\\0\\\mathcal{N}d_1
\end{smallmatrix}\right), \cdots,
\left(\begin{smallmatrix} 0 \\
\vdots\\0\\\mathcal{N}d_n
\end{smallmatrix}\right)\right)$.

\vskip5pt

{\bf Step 2.} \ Write  $Y_{(\theta_i)}=\tau_{\mathcal M} {\bf
Cok}X_{(\phi_i)}\oplus\bigoplus\limits_{i=1}^{n-1}{\bf
p}_{n-i+1}(\mathcal{N}Q'_{i})$. By taking the $i$-th branches and
the $n$-th branches of terms in $(2.7)$, we get the following
commutative diagram with exact rows
$$\xymatrix@C=0.8cm@R=0.3cm{0\ar[r]& Y_n\ar[r]^-{\left(\begin{smallmatrix}a
\\b\end{smallmatrix}\right)}\ar[d]_{\theta_i\cdots\theta_{n-1}} &
{\bigoplus\limits_{j=1}^i\mathcal{N}Q_j\oplus\bigoplus
\limits_{j=i+1}^n\mathcal{N}Q_j}\ar[d]^{(1,
0)}\ar[rrrr]^-{((\mathcal{N}d_1, \cdots, \mathcal{N}d_i), (\mathcal
Nd_{i+1},
\cdots, \mathcal{N}d_{n})) } & & & &\mathcal{N}P\ar[d] \\
0\ar[r] & Y_i\ar[r]^-{\cong} &
\bigoplus\limits_{j=1}^i\mathcal{N}Q_j\ar[rrrr] & && &0.}$$ In
particular, $Y_n=\Ker(\mathcal{N}d_{1}, \cdots, \mathcal{N}d_{n})$.
The upper exact sequence means that
$$\xymatrix@R=0.3cm@C=1cm {Y_n\ar[r]^-{-b}\ar[d]^-a & \bigoplus
\limits_{j=i+1}^n\mathcal{N}Q_j\ar[d]^-{(\mathcal{N}d_{i+1}, \cdots,
\mathcal{N}d_n)}\\
\bigoplus \limits_{j=1}^i\mathcal{N}Q_j\ar[r]^-{(\mathcal{N}d_1,
\cdots, \mathcal{N}d_i)} & \mathcal{N}P}$$ is a pull back square,
for each $1\leq i\leq n-1$. It follows that
$$\Ker (\theta_i\cdots\theta_{n-1}) = \Ker a = \Ker (\mathcal{N}d_{i+1}, \cdots,
\mathcal{N}d_n), \ 1\le i\le n-1,$$ and hence ${\bf Ker}
Y_{(\theta_i)}=\left(\begin{smallmatrix}\Ker(\mathcal{N}d_{1},
\cdots, \mathcal{N}d_{n})
\\\Ker(\mathcal{N}d_{2}, \cdots,
\mathcal{N}d_{n})\\\vdots\\\Ker(\mathcal{N}d_{n-1},
\mathcal{N}d_{n})\\\Ker\mathcal{N}d_{n}\end{smallmatrix}\right)_{(\theta'_i)}.$
We explicitly compute ${\bf Ker} Y_{(\theta_i)}$ below.

\vskip5pt

{\bf Step 3.} \ By (2.5) we get the following commutative diagram
with exact rows:
$$\xymatrix@R=0.5cm@C=0.6cm{0\ar[r] & \tau
X_1\ar[r]^-{\gamma_n}\ar[d]_-{\beta_{n-1}} &
\mathcal{N}Q_n\ar[rr]^-{\mathcal{N}d_n}\ar[d]^-{\left(\begin{smallmatrix}0\\1
\end{smallmatrix}\right)} && \mathcal{N}P\ar@{=}[d]\\
0\ar[r] &
\tau\Cok(\phi_1\cdots\phi_{n-1})\oplus\mathcal{N}Q_n\ar[r]^-{\gamma_{n-1}}\ar[d]_-{\beta_{n-2}}
&
\mathcal{N}Q_{n-1}\oplus\mathcal{N}Q_n\ar[rr]^-{(\mathcal{N}d_{n-1},\
\mathcal{N}d_n)}\ar[d]^-{\left(\begin{smallmatrix}0\\E_2
\end{smallmatrix}\right)} && \mathcal{N}P\ar@{=}[d]\\
& \vdots\ar[d] & \vdots\ar[d] && \vdots\ar@{=}[d]\\
0\ar[r] &
\tau\Cok(\phi_1\phi_2)\oplus\bigoplus\limits_{i=3}^n\mathcal{N}Q_i\ar[r]^-{\gamma_2}\ar[d]_-{\beta_{1}}
&
\mathcal{N}Q_2\oplus\bigoplus\limits_{i=3}^n\mathcal{N}Q_i\ar[rr]^-{(\mathcal{N}d_2,\cdots,
\mathcal{N}d_n)}\ar[d]^-{\left(\begin{smallmatrix}0\\E_{n-1}
\end{smallmatrix}\right)}&& \mathcal{N}P\ar@{=}[d]\\
0\ar[r] &
\tau\Cok\phi_1\oplus\bigoplus\limits_{i=2}^n\mathcal{N}Q_i\ar[r]^-{\gamma_1}
&
\mathcal{N}Q_1\oplus\bigoplus\limits_{i=2}^n\mathcal{N}Q_i\ar[rr]^-{(\mathcal{N}d_1,\cdots,
\mathcal{N}d_n)}  && \mathcal{N}P}\eqno(2.8)$$ where $E_i$ is the
identity matrix, $\gamma_i=\left(\begin{smallmatrix}
  \sigma_{i} & & (-\mathcal{N}s_{i} &
  -\mathcal{N}(s_{i}s_{i+1}) & \cdots &
  -\mathcal{N}(s_{i}\cdots s_{n-1}))\\
  0 &  & & & E_{n-i} & \\
\end{smallmatrix}\right)_{(n-i+1)\times(n-i+1)}$\\[0.6em] for $1\leq i\leq n,$ $\beta_i=\left(\begin{smallmatrix}
  \alpha_i & &(0 & 0 & \cdots & 0)\\
  \sigma_{i+1} & & (-\mathcal{N}s_{i+1} &
  -\mathcal{N}(s_{i+1}s_{i+2}) & \cdots &
  -\mathcal{N}(s_{i+1}\cdots s_{n-1}))\\
  0 &  &  & & E_{n-i-1} &\\
\end{smallmatrix}\right)_{(n-i+1)\times(n-i)}$ for $1\leq i\leq n-1.$
From (2.8) we see ${\bf
Ker}Y_{(\theta_i)}\cong\left(\begin{smallmatrix}
\tau\Cok\phi_1\oplus\bigoplus\limits_{j=2}^n\mathcal{N}Q_j\\
\tau\Cok(\phi_1\phi_2)\oplus\bigoplus\limits_{j=3}^n\mathcal{N}Q_j\\
\vdots\\
\tau\Cok(\phi_1\cdots\phi_{n-1})\oplus\mathcal{N}Q_n\\
\tau X_1
\end{smallmatrix}\right)_{(\beta_i)}.$
Applying Lemma \ref{Mimosummand}$(i)$, it is isomorphic to ${\bf
Mimo} \ \tau \ {\bf Cok}X_{(\phi_i)}\oplus J,$ where $J$ is an
injective object in $\mathcal{S}_n(A)$. Thus
$$\begin{array}{ll}
{\bf Mimo}&\tau \ {\bf Cok}X_{(\phi_i)}\oplus J \cong {\bf
Ker}Y_{(\theta_i)} \overset{_{Def.}}\cong {\bf Ker} \ \tau_{\mathcal
M} \ {\bf Cok}X_{(\phi_i)}\oplus{\bf
Ker}(\bigoplus\limits_{i=1}^{n-1}{\bf
p}_{n-i+1}(\mathcal{N}Q'_{i})) \\[1em]
& \overset{_{Cor.2.3}}\cong \tau_{\mathcal{S}}X_{(\phi_i)}\oplus{\bf
Ker}(\bigoplus\limits_{i=1}^{n-1}{\bf p}_{n-i+1}(\mathcal{N}Q'_{i}))
\overset{(1.2)}\cong \tau_{\mathcal{S}}X_{(\phi_i)}\oplus
\bigoplus\limits_{i=1}^{n-1}{\bf m}_{i}(\mathcal{N}Q'_{i}).
\end{array}$$
Since ${\bf Mimo} \ \tau \ {\bf Cok}X_{(\phi_i)}$ and
$\tau_{\mathcal{S}}X_{(\phi_i)}$ have no nonzero injective direct
summands in $\mathcal{S}_n(A)$ (cf. Remark $(ii)$ in 1.4), and
$\mathcal{S}_n(A)$ is Krull-Schmidt, we get
$\tau_{\mathcal{S}}X_{(\phi_i)} \cong {\bf Mimo} \ \tau \ {\bf
Cok}X_{(\phi_i)}.$ $\s$

\subsection{Example.} \ Let $A$,
$S$, $i$, and $\pi$ be as in 2.3. Then there are 6 indecomposable
non-projective objects in $\mathcal{S}_{3}(A)$. By Theorem 2.4 we
have
$$\begin{array}{l}\tau_{\mathcal{S}}{\left(\begin{smallmatrix}
      A\\S\\0
    \end{smallmatrix}\right)_{(0, i)}}=
    {\bf Mimo} \ \tau {\left(\begin{smallmatrix}
     S\\A\\A
    \end{smallmatrix}\right)_{(1, \pi)}}=
    {\bf Mimo}{\left(\begin{smallmatrix}
      S\\0\\0
    \end{smallmatrix}\right)_{(0, 0)}}={\left(\begin{smallmatrix}
     S\\0\\0
    \end{smallmatrix}\right)_{(0, 0)}}\\[1em]
    \tau_{\mathcal{S}}{\left(\begin{smallmatrix}
      S\\0\\0
    \end{smallmatrix}\right)_{(0, 0)}}=
    {\bf Mimo} \ \tau{\left(\begin{smallmatrix}
      S\\S\\S
    \end{smallmatrix}\right)_{(1, 1)}}=
    {\bf Mimo}{\left(\begin{smallmatrix}
      S\\S\\S
    \end{smallmatrix}\right)_{(1, 1)}}=
    {\left(\begin{smallmatrix}
      S\\S\\S
    \end{smallmatrix}\right)_{(1, 1)}}
    \end{array}$$

$$\begin{array}{l}\tau_{\mathcal{S}}{\left(\begin{smallmatrix}
      S\\S\\S
    \end{smallmatrix}\right)_{(1, 1)}}=
    {\bf Mimo}\ \tau{\left(\begin{smallmatrix}
      0\\0\\S
    \end{smallmatrix}\right)_{(0, 0)}}=
    {\bf Mimo}{\left(\begin{smallmatrix}
      0\\0\\S
    \end{smallmatrix}\right)_{(0, 0)}}=
    {\left(\begin{smallmatrix}
      A\\A\\S
    \end{smallmatrix}\right)_{(i, 1)}}\\[1em]
     \tau_{\mathcal{S}}{\left(\begin{smallmatrix}
      A\\A\\S
    \end{smallmatrix}\right)_{(i, 1)}}=
    {\bf Mimo} \ \tau{\left(\begin{smallmatrix}
      0\\S\\A
    \end{smallmatrix}\right)_{(\pi, 0)}}=
    {\bf Mimo}{\left(\begin{smallmatrix}
      0\\S\\0
    \end{smallmatrix}\right)_{(0, 0)}}=
    {\left(\begin{smallmatrix}
      A\\S\\0
    \end{smallmatrix}\right)_{(1, i)}}\\[1em]
     \tau_{\mathcal{S}}{\left(\begin{smallmatrix}
      S\\S\\0
    \end{smallmatrix}\right)_{(0, 1)}}=
    {\bf Mimo} \ \tau{\left(\begin{smallmatrix}
      0\\S\\S
    \end{smallmatrix}\right)_{(1, 0)}}=
    {\bf Mimo}{\left(\begin{smallmatrix}
      0\\S\\S
    \end{smallmatrix}\right)_{(1, 0)}}=
    {\left(\begin{smallmatrix}
      A\\S\\S
    \end{smallmatrix}\right)_{(1, i)}}\\[1em]
    \tau_{\mathcal{S}}{\left(\begin{smallmatrix}
      A\\S\\S
    \end{smallmatrix}\right)_{(1, i)}}=
    {\bf Mimo} \ \tau{\left(\begin{smallmatrix}
      S\\S\\A
    \end{smallmatrix}\right)_{(\pi, 1)}}=
    {\bf Mimo}{\left(\begin{smallmatrix}
      S\\S\\0
    \end{smallmatrix}\right)_{(0, 1)}}=
    {\left(\begin{smallmatrix}
      S\\S\\0
    \end{smallmatrix}\right)_{(0, 1)}}\end{array}$$
The Auslander-Reiten quiver of $\mathcal{S}_{3}(A)$ looks like
$$\xymatrix@C=0.6cm@R=0.05cm{& {\begin{smallmatrix}A\\0\\0\end{smallmatrix}}\ar[dr] & &
{\begin{smallmatrix}A\\A\\0\end{smallmatrix}}\ar[dr] & &
{\begin{smallmatrix}A\\A\\A\end{smallmatrix}}\ar[dr] &\\
{\begin{smallmatrix}S\\0\\0\end{smallmatrix}}\ar[ur]\ar[dr] & &
{\begin{smallmatrix}A\\S\\0\end{smallmatrix}}\ar[ur]\ar[dr]\ar@{.>}[ll]
& &
{\begin{smallmatrix}A\\A\\S\end{smallmatrix}}\ar[ur]\ar[dr]\ar@{.>}[ll]
& &
{\begin{smallmatrix}S\\S\\S\end{smallmatrix}}\ar@{.>}[ll]\\
& {\begin{smallmatrix}S\\S\\0\end{smallmatrix}}\ar[ur]\ar[dr] & &
{\begin{smallmatrix}A\\S\\S\end{smallmatrix}}\ar[ur]\ar[dr]\ar@{.>}[ll]
& &
{\begin{smallmatrix}S\\S\\0\end{smallmatrix}}\ar[ur]\ar@{.>}[ll] & \\
& & {\begin{smallmatrix}S\\S\\S\end{smallmatrix}}\ar[ur] & &
{\begin{smallmatrix}S\\0\\0\end{smallmatrix}}\ar[ur]\ar@{.>}[ll] & &
}$$

\vskip10pt

\subsection{} \ Denote by $\mathcal{S}_n(A)_{\mathcal{I}}$ the full subcategory of
$\mathcal{S}_n(A)$ consisting of all the objects which have no
nonzero injective direct summands of $\mathcal{S}_n(A)$. By Theorem
\ref{tauformula}, we have

\begin{cor}\label{Mimoform}
\ Every object in $\mathcal{S}_n(A)_{\mathcal{I}}$ has the form
${\bf Mimo}X,$ where each $X_i$ has no nonzero injective direct
summands.
\end{cor}

The following result will be used in the next section, whose proof
is omitted, since it is the same as the case of $n=2$ (see [RS2],
Corollary 5.4).

\vskip10pt
\begin{cor}\label{functorF}
\ The canonical functor $W:
\mathcal{S}_n(A)_{\mathcal{I}}\rightarrow {\rm
Mor}_n(A\mbox{-}\overline {\rm mod})$ given by $X_{(\phi_i)}\mapsto
X_{(\overline{\phi}_i)}$ is dense, preserves indecomposables, and
reflects isomorphisms.
\end{cor}

\subsection{} Dual to Theorem 2.4 we have

\begin{thm}  Let $X\in \mathcal{F}_n(A)$. Then $\tau_{\mathcal{F}} X \cong {\bf Cok} \ {\bf Mimo} \ \tau X,$
and $\tau_{\mathcal{F}}^-X \cong {\bf Mepi} \ \tau^- \ {\bf Ker}X.$
\end{thm}

\section{\bf Applications to selfinjective algebras}

Throughout this section, $A$ is a selfinjective algebra. Then
$A\mbox{-}\overline{\rm mod} = A\mbox{-}\underline{\rm mod}$ is a
triangulated category with the suspension functor $\Omega^{-1}$
([H], p.16), where $\Omega^{-1}$ is the cosyzygy of $A$. The
distinguished triangles of $A\mbox{-}\underline{\rm mod}$ are
exactly triangles isomorphic to those given by all the short exact
sequences in $A$-mod. Note that $\Omega$ and $\mathcal N$ commute
and $\tau \cong \Omega^2 \mathcal{N}\cong \mathcal{N}\Omega^2$ is an
endo-equivalence of $A$-$\underline{\rm mod}$  (\cite{ARS}, p.126),
and that $\tau$ is a triangle functor.

A rotation of an object in ${\rm Mor}_2(A\mbox{-}\underline {\rm
mod})$ is introduced by Ringel and Schmidmeier [RS2]. The definition
of a rotation of an object in ${\rm Mor}_n(A\mbox{-}\underline {\rm
mod})$ needs new considerations. We have to take up pages to justify
that it is well-defined. Then we get a formula for $\underline
{\tau_{\mathcal S}^jX}\in {\rm Mor}_n(A\mbox{-}\underline {\rm
mod})$ for $X\in \mathcal S_n(A)$ and $j\ge 1$. This is applied to
the study of the periodicity of $\tau_{\mathcal S}$ on the objects
of $\mathcal S_n(A)$. In particular, for the selfinjective Nakayama
algebras $\A(m, t)$ we have $\tau_{\mathcal S}^{2m(n+1)}X\cong X$
for $X\in \mathcal{S}_n(\A(m, t))$.

\subsection{} \  Let $X_{(\underline {\phi_i})}\in {\rm Mor}_n(A\mbox{-}\underline
{\rm mod})$. Just choose $\phi_i$ as representatives for the
morphisms $\underline {\phi_i}$ in  $A\mbox{-}\underline {\rm mod}$.
Let $h_{i+1}: X_{i+1}\rightarrow I_{i+1}$ be an injective envelope
with cokernel $\Omega^{-1}X_{i+1}$, $1\leq i\leq n-1$. Taking
pushout we get the following commutative diagram with exact rows
$$\xymatrix@R=0.4cm@C=1cm{0\ar[r] & X_{i+1}\ar[r]^-{h_{i+1}}\ar[d]_{\phi_1\cdots\phi_{i}}
& I_{i+1}\ar[r]\ar[d]^{j_{i+1}} & \Omega^{-1}X_{i+1}\ar[r]\ar@{=}[d] & 0\\
0\ar[r]& X_1\ar[r]^-{g_{i+1}} & Y^1_{i+1}\ar[r] &
\Omega^{-1}X_{i+1}\ar[r] &  0.}$$ This gives the exact sequence
$$\xymatrix@1{0\ar[r] & X_{i+1}\ar[rr]^-{\left(\begin{smallmatrix}
  \phi_1\cdots\phi_{i} \\ h_{i+1}
\end{smallmatrix}\right)}& &X_1\oplus I_{i+1}\ar[rr]^-{(g_{i+1},\
-j_{i+1})}& & Y^1_{i+1}\ar[r] & 0,}\eqno(3.1)$$ and induces the
following commutative diagram with exact rows, $1\leq i\leq n-2$
$$\xymatrix@R=0.4cm@C=1cm{0\ar[r] & X_{i+2}\ar[rr]^-{\left(\begin{smallmatrix}
  \phi_1\cdots\phi_{i+1} \\ h_{i+2}
\end{smallmatrix}\right)}\ar[d]_-{\phi_{i+1}}& & X_1\oplus I_{i+2}\ar[rr]^-{(g_{i+2},\
-j_{i+2})}\ar@{.>}[d]^{\left(\begin{smallmatrix}
  1&0 \\ 0 & *
\end{smallmatrix}\right)} & & Y^1_{i+2}\ar[r]\ar@{.>}[d]^-{\psi_{i}}
& 0\\0\ar[r] & X_{i+1}\ar[rr]^-{\left(\begin{smallmatrix}
  \phi_1\cdots\phi_{i} \\ h_{i+1}
\end{smallmatrix}\right)}& &X_1\oplus I_{i+1}\ar[rr]^-{(g_{i+1},\
-j_{i+1})}& & Y^1_{i+1}\ar[r] & 0.}\eqno(3.2)$$ By (3.2) we have
$g_{i+1}=\psi_ig_{i+2}, \ 1\le i\le n-2$. Put $\psi_{n-1}=g_n$. Then
$g_{i+1}=\psi_i\cdots\psi_{n-1}, \  1\leq i\leq n-1$. By the
construction of a distinguished triangle in $A\mbox{-}\underline{\rm
mod}$, we get distinguished triangles from (3.1)
$$X_{i+1} \stackrel{\underline{\phi_1}\cdots \underline{\phi_i}} \longrightarrow X_1\rightarrow Y_{i+1}^1\rightarrow
\Omega^{-1}X_{i+1}, \ 1\leq i \leq n-1, \eqno(3.3)$$ and by (3.2) we
get the following commutative diagram, where the rows are
distinguished triangles from (3.3)
$$\xymatrix@R=0.4cm@C=0.9cm{X_n\ar[r]\ar[d]^-{\underline{\phi_{n-1}}} & X_1\ar@{=}[d]\ar[r]^-{\underline{\psi_{n-1}}} &
Y_n^1\ar[r]\ar@{.>}[d]^-{\underline{\psi_{n-2}}} & \Omega^{-1}X_n\ar[d]\\
X_{n-1}\ar[r]\ar[d]^-{\underline{\phi_{n-2}}} & X_1\ar@{=}[d]\ar[r]
& Y_{n-1}^1\ar[r]\ar@{.>}[d]^-{\underline{\psi_{n-3}}} &
\Omega^{-1}X_{n-1}\ar[d]\\
\vdots\ar[d]^-{\underline{\phi_2}} & \vdots
\ar@{=}[d] & \vdots\ar@{.>}[d]^-{\underline{\psi_1}} & \vdots\ar[d]\\
X_2\ar[r]^{\underline{\phi_1}} & X_1\ar[r] & Y_2^1\ar[r] &
\Omega^{-1}X_2.}\eqno(3.4)$$ The {\it rotation} ${\bf Rot}
X_{(\underline {\phi_i})}$ of $X_{(\underline{\phi_i})}$ is defined
to be
$$\xymatrix@1{(X_1\ar[r]^{\underline{\psi_{n-1}}}&Y_n^1\ar[r]&
\cdots\ar[r]^{\underline {\psi_1}}&Y_2^1)}\in {\rm
Mor}_n(A\mbox{-}\underline {\rm mod})$$ (here and in the following,
for convenience we write the rotation in a row). We remark  that
${\bf Rot} X_{(\underline {\phi_i})}$ is well-defined: if
$X_{(\underline{\phi_i})}\cong Y_{(\underline{\theta_i})}$ in ${\rm
Mor}_n(A\mbox{-}\underline {\rm mod})$ with all $X_i$ and $Y_i$
having no nonzero injective direct summands, then ${\bf
Mimo}X_{({\phi_i})}\cong{\bf Mimo}Y_{({\theta_i})}$ by Lemma
1.6$(i)$, and hence ${\bf Rot}X_{(\underline{\phi_i})} \cong {\bf
Rot}Y_{(\underline{\theta_i})}$, by Lemma 3.1 below.

\vskip10pt

\begin{lem}\label{rotationlemma}
\ Let $X_{(\phi_i)}\in {\rm Mor}_n(A)$. Then ${\bf Rot}
X_{(\underline {\phi_i})}\cong\underline{{\bf Cok \
Mimo}X_{(\phi_i)}}$ in ${\rm Mor}_n(A\mbox{-}\underline {\rm mod})$.
\end{lem}

\vskip10pt

Before proving Lemma 3.1, for later convenience, we restate Claim 2
in $\S4$ of [RS2] in the more explicit way we will use.

\vskip10pt

\begin{lem} \ Let $\xymatrix@=0.5cm{0\ar[r] & A\ar[r]^-{\left(\begin{smallmatrix}f\\
h \end{smallmatrix}\right)} & B\oplus I\ar[r] & C\ar[r] & 0}$ be an
exact sequence with $I$ an injective $A$-module. Then there is an
injective $A$-module $J$ such that $I={\rm IKer}f\oplus J$, and that
the following diagram with exact rows commutes
$$\xymatrix{0\ar[r] & A\ar@{=}[d]\ar[r]^-{\left(\begin{smallmatrix}f\\
e\\0\end{smallmatrix}\right)} & B\oplus{\rm IKer}f\oplus
J\ar[d]^-{\left(\begin{smallmatrix}1 & 0& 0\\ 0 & 1 & 0
\\w & 0 & 1\end{smallmatrix}\right)}\ar[r] & C'\oplus J\ar[d]^-{\wr}\ar[r] & 0\\
0\ar[r] & A\ar[r]^-{\left(\begin{smallmatrix}f\\
e\\ h' \end{smallmatrix}\right)} & B\oplus{\rm IKer}f\oplus J\ar[r]
& C\ar[r] & 0,}$$ where $h=\left(\begin{smallmatrix} e\\ h'
\end{smallmatrix}\right)$, $e: A\rightarrow{\rm
IKer}f$ is an extension of the injective envelope $\Ker
f\hookrightarrow{\rm IKer}f$, $h': A\rightarrow J$ satisfies $h'\Ker
f=0$,  and $C'=\Cok\left(\begin{smallmatrix}f\\
e\end{smallmatrix}\right)$.
\end{lem}

\vskip10pt

{\bf Proof of Lemma 3.1.} We divide the proof into three steps.

{\bf Step 1.} \ Recall ${\bf
Mimo}X_{(\phi_i)}=\left(\begin{smallmatrix}
X_1\oplus\bigoplus\limits_{l=1}^{n-1}{\rm IKer}\phi_l\\
\vdots \\X_{n-1}\oplus\ik\phi_{n-1}\\X_n
\end{smallmatrix}\right)_{(\theta_i)}  {\rm with} \
 \theta_i = \left(\begin{smallmatrix}
\phi_i& 0 & \cdots &  0 \\
e_i& 0 & \cdots &  0 \\
0 & 1& \cdots &  0 \\
\vdots & \vdots &
\cdots & \vdots \\
0 & 0& \cdots & 1
\end{smallmatrix}\right)_{(n-i+1)\times (n-i)},$
$e_i: X_{i+1}\rightarrow \ik\phi_i$ is an extension of the injective
envelope $\Ker\phi_i\hookrightarrow \ik\phi_i$,  and $${\bf Cok \
Mimo}X_{(\phi_i)}=\left(\begin{smallmatrix}
\Cok(\theta_1)\\ \Cok(\theta_1\theta_2) \\ \vdots\\ \Cok(\theta_1\cdots\theta_{n-1}) \\
X_1\oplus\bigoplus\limits_{l=1}^{n-1}{\rm IKer}\phi_l
\end{smallmatrix}\right)_{(\theta'_i)}.$$ Since
$\theta_1\cdots \theta_i = diag(\alpha_{i}, E_{n-i-1}):
X_{i+1}\oplus\bigoplus\limits_{l=i+1}^{n-1}{\rm IKer}\phi_l
\rightarrow X_{1}\oplus\bigoplus\limits_{l=1}^{n-1}{\rm
IKer}\phi_l,$
where $\alpha_{i} = \left (\begin{smallmatrix}\phi_1\cdots\phi_i\\e_1\phi_2\cdots\phi_i\\
\vdots\\e_{i-1}\phi_i\\e_i\end{smallmatrix}\right): X_{i+1}
\rightarrow X_{1}\oplus\bigoplus\limits_{l=1}^{i}{\rm IKer}\phi_l$,
and $E_{n-i-1}$ is the identity matrix, we get the following
commutative diagram with exact rows, $2\leq i\leq n-1$,
$$\xymatrix@R=0.4cm@C=0.8cm{0\ar[r] &  X_{i+1}\ar[r]^-{\alpha_i}\ar[d]_-{\phi_i} &
X_1\oplus\bigoplus\limits_{l=1}^{i}{\rm
IKer}\phi_l\ar[r]^-{\pi_i}\ar[d]^-{(E_i,0)}
&\Cok(\theta_1\cdots\theta_i)\ar[r]\ar[d]^-{\theta'_{i-1}} & 0\\
0\ar[r]& X_i\ar[r]^-{\alpha_{i-1}}&
X_1\oplus\bigoplus\limits_{l=1}^{i-1}{\rm
IKer}\phi_l\ar[r]^-{\pi_{i-1}}&
\Cok(\theta_1\cdots\theta_{i-1})\ar[r] & 0}\eqno(3.5)$$ with
$\pi_{n-1} = \theta'_{n-1}$.

Applying  Lemma 3.2 to the upper exact sequence of (3.5) for $1\leq
i\leq n-1$, we get injective $A$-modules $J_{i+1}$ such that
$\bigoplus\limits_{l=1}^{i}{\rm IKer}\phi_l = {\rm
IKer}(\phi_1\cdots\phi_i)\oplus J_{i+1}$ and that the following
diagram with exact rows commutes
$$\xymatrix@C=0.9cm@R=0.7cm{0\ar[r]&
X_{i+1}\ar@{=}[d]\ar[r]^-{\left(\begin{smallmatrix}\phi_1\cdots\phi_i\\a_i\\0\end{smallmatrix}\right)}
& X_1\oplus{\rm IKer}(\phi_1\cdots\phi_i)\oplus J_{i+1}
\ar[d]^-{\left(\begin{smallmatrix}1&0&0\\0&1&0\\w_i&0&1\end{smallmatrix}\right)}
\ar[r]^-{\left(\begin{smallmatrix}b_i^1&b_i^2&0\\0&0&1\end{smallmatrix}\right)}
& Z_{i+1}\oplus J_{i+1}\ar[r]\ar[d]_-{\wr}^-{\beta_i} & 0\\
0\ar[r]&
X_{i+1}\ar[r]^-{\left(\begin{smallmatrix}\phi_1\cdots\phi_i\\a_i\\d_i\end{smallmatrix}\right)}&
X_1\oplus{\rm IKer}(\phi_1\cdots\phi_i)\oplus J_{i+1}\ar[r]^-{\pi_i}
& \Cok(\theta_1\cdots\theta_i)\ar[r] & 0}\eqno(3.6)$$ where
$\alpha_i=\left(\begin{smallmatrix}\phi_1\cdots\phi_i\\a_i\\d_i\end{smallmatrix}\right)$,
$a_i: X_{i+1}\rightarrow{\rm IKer}(\phi_1\cdots\phi_i)$ is an
extension of the injective envelope
$\Ker(\phi_1\cdots\phi_i)\hookrightarrow{\rm
IKer}(\phi_1\cdots\phi_i)$, $d_i: X_{i+1}\rightarrow J_{i+1}$
satisfies $d_i\Ker(\phi_1\cdots\phi_i)$=0, and $Z_{i+1} = \Cok
\left(\begin{smallmatrix}\phi_1\cdots\phi_i\\a_i\end{smallmatrix}\right)$.
Thus by (3.6) and (3.5) we get the following commutative diagram
with exact rows for $2\le i\le n-1$ (where the two rows in the
middle come from (3.5)):
$$\xymatrix@C=0.6cm@=0.7cm{0\ar[r]&
X_{i+1}\ar@{=}[d]\ar[rr]^-{\left(\begin{smallmatrix}\phi_1\cdots\phi_i\\a_i\\0\end{smallmatrix}\right)}
&& X_1\oplus{\rm IKer}(\phi_1\cdots\phi_i)\oplus J_{i+1}
\ar[d]^-{\left(\begin{smallmatrix}1&0&0\\ 0 & 1 & 0\\ w_i & 0 & 1
\end{smallmatrix}\right)}
\ar[rr]^-{\left(\begin{smallmatrix}b_i^1&b_i^2&0\\0&0&1\end{smallmatrix}\right)}
&& Z_{i+1}\oplus J_{i+1}\ar[r]\ar[d]^-{\beta_i}_-{\wr}
&0\\
0\ar[r]&X_{i+1}\ar[rr]^-{\left(\begin{smallmatrix}\phi_1\cdots\phi_i\\a_i\\d_i\end{smallmatrix}\right)}\ar[d]_-{\phi_i}&&
X_1\oplus{\rm IKer}(\phi_1\cdots\phi_i)\oplus
J_{i+1}\ar[rr]^-{\pi_i}\ar[d]^-{\left(\begin{smallmatrix}1&0&0\\ 0 & * & *\\
0 & * & *
\end{smallmatrix}\right)} &&
\Cok(\theta_1\cdots\theta_i)\ar[r]\ar[d]^-{\theta'_{i-1}} & 0\\
0\ar[r]&X_i\ar[rr]^-{\left(\begin{smallmatrix}\phi_1\cdots\phi_{i-1}\\a_{i-1}\\d_{i-1}\end{smallmatrix}\right)}\ar@{=}[d]&&
X_1\oplus{\rm IKer}(\phi_1\cdots\phi_{i-1})\oplus
J_i\ar[rr]^-{\pi_{i-1}}\ar[d]^-{\left(\begin{smallmatrix}1&0&0\\ 0 & 1 & 0\\
w_{i-1} & 0 & 1
\end{smallmatrix}\right)^{-1}} &&
\Cok(\theta_1\cdots\theta_{i-1})\ar[r]\ar[d]^-{\beta_{i-1}^{-1}}_-{\wr} & 0\\
0\ar[r]&X_i\ar[rr]^-{\left(\begin{smallmatrix}\phi_1\cdots\phi_{i-1}\\a_{i-1}\\0\end{smallmatrix}\right)}
&& X_1\oplus{\rm IKer}(\phi_1\cdots\phi_{i-1})\oplus J_i
\ar[rr]^-{\left(\begin{smallmatrix}b_{i-1}^1&b_{i-1}^2&0\\0&0&1\end{smallmatrix}\right)}
&& Z_i\oplus J_i\ar[r] & 0.}$$  Taking the first and the last rows,
we get the following commutative diagram with exact rows:
$$\xymatrix@C=0.7cm@R=0.7cm{0\ar[r]&
X_{i+1}\ar[d]_-{\phi_i}\ar[rr]^-{\left(\begin{smallmatrix}\phi_1\cdots\phi_i\\a_i\\0\end{smallmatrix}\right)}
&& X_1\oplus{\rm IKer}(\phi_1\cdots\phi_i)\oplus J_{i+1}
\ar[d]^-{\left(\begin{smallmatrix}1&0&0\\c_{i1}&c_{i2}&
*\\ *& * & *\end{smallmatrix}\right)}
\ar[rr]^-{\left(\begin{smallmatrix}b_i^1&b_i^2&0\\0&0&1\end{smallmatrix}\right)}
&& Z_{i+1}\oplus
J_{i+1}\ar[r]\ar[d]^-{\beta^{-1}_{i-1}\theta'_{i-1}\beta_i} & 0\\
0\ar[r]&
X_i\ar[rr]^-{\left(\begin{smallmatrix}\phi_1\cdots\phi_{i-1}\\a_{i-1}\\0\end{smallmatrix}\right)}
&& X_1\oplus{\rm IKer}(\phi_1\cdots\phi_{i-1})\oplus J_i
\ar[rr]^-{\left(\begin{smallmatrix}b_{i-1}^1&b_{i-1}^2&0\\0&0&1\end{smallmatrix}\right)}
&& Z_i\oplus J_i\ar[r] & 0,}$$ where for later convenience we write
$$\beta^{-1}_{i-1}\theta'_{i-1}\beta_i =
\left(\begin{smallmatrix}f_{i-1}&*\\
*&*\end{smallmatrix}\right), \ 2\le i\le n-1.\eqno(3.7)$$
Taking off $J_i$ and $J_{i+1}$, we get the following commutative
diagram with exact rows, $2\leq i\leq n-1$,
$$\xymatrix@C=0.65cm@R=0.8cm{0\ar[r]&
X_{i+1}\ar[d]_-{\phi_i}\ar[rr]^-{\left(\begin{smallmatrix}\phi_1\cdots\phi_i\\a_i\end{smallmatrix}\right)}
&& X_1\oplus{\rm IKer}(\phi_1\cdots\phi_i)
\ar[d]^-{\left(\begin{smallmatrix}1&0\\c_{i1}&c_{i2}\end{smallmatrix}\right)}
\ar[rr]^-{(b_i^1, \ b_i^2)}
&& Z_{i+1}\ar[r]\ar[d]^-{f_{i-1}} & 0\\
0\ar[r]&
X_i\ar[rr]^-{\left(\begin{smallmatrix}\phi_1\cdots\phi_{i-1}\\a_{i-1}\end{smallmatrix}\right)}
&& X_1\oplus{\rm IKer}(\phi_1\cdots\phi_{i-1}) \ar[rr]^-{(b_{i-1}^1,
\ b_{i-1}^2)} && Z_i\ar[r] & 0.}\eqno(3.8)$$

{\bf Step 2.}  \ Now we consider the rotation ${\bf Rot}
X_{(\underline {\phi_i})}$. Recall from the beginning of this
subsection that $h_{i+1}: X_{i+1}\rightarrow I_{i+1}$ is an
injective envelope.  For $1\leq i\leq n-1$, applying Lemma 3.2 to
(3.1), we get injective $A$-modules $J'_{i+1}$ such that
$I_{i+1}={\rm IKer}(\phi_1\cdots\phi_i)\oplus J'_{i+1}$ and that the
following diagram with exact rows commutes (cf. (3.6))
$$\xymatrix@C=0.6cm{0\ar[r]&
X_{i+1}\ar[rr]^-{\left(\begin{smallmatrix}\phi_1\cdots\phi_i\\a_i\\0\end{smallmatrix}\right)}
\ar@{=}[d] & &X_1\oplus {\rm IKer}(\phi_1\cdots\phi_i)\oplus
J'_{i+1}\ar[rr]^-{\left(\begin{smallmatrix}b_i^1&b_i^2&0\\0&0&1\end{smallmatrix}\right)}
\ar[d]^-{\left(\begin{smallmatrix}1&0&0\\0&1&0\\w'_i&0&1\end{smallmatrix}\right)}
& &
Z_{i+1}\oplus J'_{i+1}\ar[r]\ar[d]_{\wr}^-{\beta'_{i}} &0\\
0\ar[r] &
X_{i+1}\ar[rr]^-{\left(\begin{smallmatrix}\phi_1\cdots\phi_i\\a_i\\d'_i\end{smallmatrix}\right)}
& & X_1\oplus {\rm IKer}(\phi_1\cdots\phi_i)\oplus
J'_{i+1}\ar[rr]^-{(g_{i+1},\ -j_{i+1})} && Y^1_{i+1}\ar[r]&
0,}\eqno(3.9)$$ where
$h_{i+1}=\left(\begin{smallmatrix}a_i\\d'_i\end{smallmatrix}\right)$,
and $d'_i: X_{i+1}\rightarrow J'_{i+1}$, satisfying
$d'_i\Ker(\phi_1\cdots\phi_i)=0$. Thus by (3.9) and (3.2) we get the
following commutative diagram with exact rows for $2\leq i\leq n-1$
(where the two rows in the middle come from (3.2)):
$$\xymatrix@C=0.7cm@R=0.7cm{0\ar[r]&
X_{i+1}\ar@{=}[d]\ar[rr]^-{\left(\begin{smallmatrix}\phi_1\cdots\phi_i\\a_i\\0\end{smallmatrix}\right)}
&& X_1\oplus{\rm IKer}(\phi_1\cdots\phi_i)\oplus J'_{i+1}
\ar[d]^-{\left(\begin{smallmatrix}1&0&0\\ 0 & 1 & 0\\ w'_i & 0 & 1
\end{smallmatrix}\right)}
\ar[rr]^-{\left(\begin{smallmatrix}b_i^1&b_i^2&0\\0&0&1\end{smallmatrix}\right)}
&& Z_{i+1}\oplus J'_{i+1}\ar[r]\ar[d]^-{\beta'_i}_-{\wr}
&0\\
0\ar[r]&X_{i+1}\ar[rr]^-{\left(\begin{smallmatrix}\phi_1\cdots\phi_i\\a_i\\d'_i\end{smallmatrix}\right)}\ar[d]_-{\phi_i}&&
X_1\oplus{\rm IKer}(\phi_1\cdots\phi_i)\oplus
J'_{i+1}\ar[rr]^-{(g_{i+1},\
-j_{i+1})}\ar[d]^-{\left(\begin{smallmatrix}1&0&0\\ 0 & * & *\\
0 & * & *
\end{smallmatrix}\right)} &&
Y_{i+1}^1\ar[r]\ar[d]^-{\psi_{i-1}} & 0\\
0\ar[r]&X_i\ar[rr]^-{\left(\begin{smallmatrix}\phi_1\cdots\phi_{i-1}\\a_{i-1}\\d'_{i-1}\end{smallmatrix}\right)}\ar@{=}[d]&&
X_1\oplus{\rm IKer}(\phi_1\cdots\phi_{i-1})\oplus
J'_i\ar[rr]^-{(g_{i},\
-j_{i})}\ar[d]^-{\left(\begin{smallmatrix}1&0&0\\ 0 & 1 & 0\\
w'_{i-1} & 0 & 1
\end{smallmatrix}\right)^{-1}} &&
Y_i^1\ar[r]\ar[d]^-{\beta'^{-1}_{i-1}}_-{\wr} & 0\\
0\ar[r]&X_i\ar[rr]^-{\left(\begin{smallmatrix}\phi_1\cdots\phi_{i-1}\\a_{i-1}\\0\end{smallmatrix}\right)}
&& X_1\oplus{\rm IKer}(\phi_1\cdots\phi_{i-1})\oplus J'_i
\ar[rr]^-{\left(\begin{smallmatrix}b_{i-1}^1&b_{i-1}^2&0\\0&0&1\end{smallmatrix}\right)}
&& Z_i\oplus J'_i\ar[r] & 0.}$$ Taking the first and the last rows,
we get the following commutative diagram with exact rows:\\
$$\xymatrix@C=0.7cm@R=0.8cm{0\ar[r]&
X_{i+1}\ar[d]_-{\phi_i}\ar[rr]^-{\left(\begin{smallmatrix}\phi_1\cdots\phi_i\\a_i\\0\end{smallmatrix}\right)}
&& X_1\oplus{\rm IKer}(\phi_1\cdots\phi_i)\oplus J'_{i+1}
\ar[d]^-{\left(\begin{smallmatrix}1&0&0\\c'_{i1}&c'_{i2}&
*\\ *& * & *\end{smallmatrix}\right)}
\ar[rr]^-{\left(\begin{smallmatrix}b_i^1&b_i^2&0\\0&0&1\end{smallmatrix}\right)}
&& Z_{i+1}\oplus
J'_{i+1}\ar[r]\ar[d]^-{\beta'^{-1}_{i-1}\psi_{i-1}\beta'_i} & 0\\
0\ar[r]&
X_i\ar[rr]^-{\left(\begin{smallmatrix}\phi_1\cdots\phi_{i-1}\\a_{i-1}\\0\end{smallmatrix}\right)}
&& X_1\oplus{\rm IKer}(\phi_1\cdots\phi_{i-1})\oplus J'_i
\ar[rr]^-{\left(\begin{smallmatrix}b_{i-1}^1&b_{i-1}^2&0\\0&0&1\end{smallmatrix}\right)}
&& Z_i\oplus J'_i\ar[r] & 0,}$$ where for later convenience we write
$$\beta'^{-1}_{i-1}\psi'_{i-1}\beta'_i =
\left(\begin{smallmatrix}f'_{i-1}&*\\
*&*\end{smallmatrix}\right), \ 2\le i\le n-1.\eqno(3.10)$$ Taking off $J'_i$ and $J'_{i+1}$, we
get the following commutative diagram with exact rows for $2\leq
i\leq n-1$:
$$\xymatrix@C=0.7cm@R=0.7cm{0\ar[r]&
X_{i+1}\ar[d]_-{\phi_i}\ar[rr]^-{\left(\begin{smallmatrix}\phi_1\cdots\phi_i\\a_i\end{smallmatrix}\right)}
&& X_1\oplus{\rm IKer}(\phi_1\cdots\phi_i)
\ar[d]^-{\left(\begin{smallmatrix}1&0\\c'_{i1}&c'_{i2}\end{smallmatrix}\right)}
\ar[rr]^-{(b_i^1, \ b_i^2)}
&& Z_{i+1}\ar[r]\ar[d]^-{f'_{i-1}} & 0\\
0\ar[r]&
X_i\ar[rr]^-{\left(\begin{smallmatrix}\phi_1\cdots\phi_{i-1}\\a_{i-1}\end{smallmatrix}\right)}
&& X_1\oplus{\rm IKer}(\phi_1\cdots\phi_{i-1}) \ar[rr]^-{(b_{i-1}^1,
\ b_{i-1}^2)} && Z_i\ar[r] & 0.}$$ Comparing the above diagram with
$(3.8)$, by a computation we easily see that $f_i-f'_i$ factors
through an injective $A$-module for each $1\leq i\leq n-2$.

\vskip10pt

{\bf Step 3.} \ Now we get the following diagram, where the first
row can be considered as ${\bf Cok \ Mimo} X_{(\phi_i)}$ (we
identify $X_1\oplus\bigoplus\limits_{l=1}^{n-1}{\rm IKer}\phi_l$
with $X_1\oplus{\rm IKer}(\phi_1\cdots\phi_{n-1})\oplus J_n$; and
identify $X_1\oplus I_n$ with $X_1\oplus {\rm
IKer}(\phi_1\cdots\phi_{n-1})\oplus J'_n$):
$$\xymatrix@C=0.4cm@R=0.7cm{X_1\oplus\bigoplus\limits_{l=1}^{n-1}{\rm
IKer}\phi_l\ar[rr]^-{\theta'_{n-1}}\ar[d]_-{\wr}^-{\left(\begin{smallmatrix}1&0&0\\
0 & 1 & 0\\ w_{n-1} & 0 & 1
\end{smallmatrix}\right)^{-1}} &&
\Cok(\theta_1\cdots\theta_{n-1})\ar[r]^-{\theta'_{n-2}}\ar[d]^-{\beta_{n-1}^{-1}}_-{\wr}
& \cdots\ar[r] &
\Cok(\theta_1\theta_2)\ar[r]^-{\theta'_1}\ar[d]^-{\beta_2^{-1}}_-{\wr}
&
\Cok\theta_1\ar[d]^-{\beta_1^{-1}}  \\
X_1\oplus{\rm IKer}(\phi_1\cdots\phi_{n-1})\oplus
J_n\ar[rr]^-{\left(\begin{smallmatrix}b_{n-1}^1&b_{n-1}^2&0\\0&0&1\end{smallmatrix}\right)}
\ar[d]^-{\left(\begin{smallmatrix}1& 0 & 0 \\
0&1&0\\ 0&0&0\end{smallmatrix}\right)}
&& Z_n\oplus J_n\ar[r]^-{\left(\begin{smallmatrix}f_{n-2}& * \\
*&*\end{smallmatrix}\right)}\ar[d]^-{\left(\begin{smallmatrix}1& 0 \\ 0&0\end{smallmatrix}\right)} & \cdots\ar[r] &
Z_3\oplus J_3\ar[r]^-{\left(\begin{smallmatrix}f_1& * \\
*&*\end{smallmatrix}\right)} \ar[d]_-{\left(\begin{smallmatrix}1& 0
\\ 0&0\end{smallmatrix}\right)} &
 Z_2\oplus J_2\ar[d]_-{\left(\begin{smallmatrix}1& 0 \\ 0&0\end{smallmatrix}\right)}\\
X_1\oplus{\rm IKer}(\phi_1\cdots\phi_{n-1})\oplus
J'_n\ar[rr]_-{\left(\begin{smallmatrix}b_{n-1}^1&b_{n-1}^2&0\\0&0&1\end{smallmatrix}\right)}\ar[d]_-{\wr}^-{\left(\begin{smallmatrix}1&0&0\\
0 & 1 & 0\\ w'_{n-1} & 0 & 1
\end{smallmatrix}\right)}
&& Z_n\oplus J'_n\ar[r]_-{\left(\begin{smallmatrix}f'_{n-2}& * \\
*&*\end{smallmatrix}\right)}\ar[d]^-{\beta'_{n-1}}_-{\wr}  & \cdots\ar[r] &
Z_3\oplus J'_3\ar[r]_-{\left(\begin{smallmatrix}f'_1& * \\
*&*\end{smallmatrix}\right)}\ar[d]^-{\beta'_2}_-{\wr} &
Z_2\oplus J'_2\ar[d]^-{\beta'_1}_-{\wr}\\
X_1\oplus I_n\ar[rr]^-{(\psi_{n-1}, -j_n)}& &
Y_n^1\ar[r]^-{\psi_{n-2}} & \cdots\ar[r] & Y_3^1\ar[r]^-{\psi_1} &
Y_2^1}$$ (note that the squares in the first two rows commute in
$A$-mod: the left square comes from (3.6); and the remaining
commutative squares come from (3.7). Also, note that the squares in
the last two rows commute in $A$-mod: the left square comes from
(3.9); and the remaining commutative squares come from (3.10)).
However, the squares in the middle may {\bf not} commute in $A$-mod;
and the point is that they commute in $A$-\underline{mod}, as we
explain below.

Note that the left square in the middle commutes by a direct
computation.
 Since $f_i-f'_i$ factors through an
injective $A$-module, $1\leq i\leq n-2$, we realize that the
remaining $n-2$ squares in the middle commute in $A$-{\underline{\rm
mod}}. It follows that the above diagram commutes in
$A$-{\underline{\rm mod}}. It is clear that the vertical morphisms
are isomorphisms in $A$-{\underline{\rm mod}}. Regarding the above
diagram in $A\mbox{-}\underline {\rm mod}$, the first row is exactly
$\underline{{\bf Cok \ Mimo}X_{(\phi_i)}}$, and the last row \\[0.6em]is
exactly ${\bf Rot} X_{(\underline{\phi_i})}$. Thus, ${\bf Rot}
X_{(\underline{\phi_i})}\cong\underline{{\bf Cok \
Mimo}X_{(\phi_i)}}$ in ${\rm Mor}_n(A\mbox{-}\underline {\rm mod})$.
This completes the proof. $\s$

\vskip 10pt

\subsection{} Let  $X_{(\phi_i)}\in
{\rm Mor}_n(A)$.  For $1\leq k <i<j\leq n$, by (3.3) and the
Octahedral Axiom we get the following commutative diagram with first
two rows and the last two columns being distinguished triangles in
$A$-$\underline{\rm mod}$:
\[\xymatrix@R=0.4cm@C=1cm{\Omega Y_j^i\ar[r]\ar[d] & X_j\ar@{=}[d]\ar[r]^{\underline
{\phi_i}\cdots\underline {{\phi}_{j-1}}} &
X_i\ar[r]\ar[d]^{ \underline {\phi_k} \cdots \underline {\phi_{i-1}}} & Y_j^i\ar[d]\\
\Omega Y_j^k\ar[r] &X_j\ar[r]^{\underline{\phi_k} \cdots\underline {\phi_{j-1}}} & X_k\ar[r]\ar[d] & Y_j^k\ar[d]\\
& & Y_i^k\ar[d]\ar@{=}[r] & Y_i^k\ar[d] & \\ && \Omega^{-1}X_i\ar[r]
& \Omega^{-1}Y_j^i.}\eqno(3.11)\] For $1\leq m\leq n$ we prove the
following formula by induction $$\begin{array}{ll} {\bf Rot}^m
X_{(\underline{\phi_i})}= & (\Omega^{-(m-2)}Y_m^{m-1}\rightarrow
\Omega^{-(m-2)}Y_m^{m-2}\rightarrow \cdots\rightarrow
\Omega^{-(m-2)}Y_m^1 \rightarrow  \Omega^{-(m-1)}X_m \rightarrow
\\ & \ \Omega^{-(m-1)}Y_n^m \ \rightarrow \
\Omega^{-(m-1)}Y_{n-1}^m \ \rightarrow \cdots\rightarrow
\Omega^{-(m-1)}Y_{m+1}^m) \ \ \ \ \ \ \ \ \ \ \ \ \ \ \ (3.12)
\end{array} $$  Convention about $(3.12)$: $\Omega^{-(m-1)}X_m$ is the $(n-m+1)$-st
branch of ${\bf Rot}^m X_{(\underline{\phi_i})}$.

By definition (3.12) holds for $m=1$. Assume that it holds for $1\le
m\le n-1.$ Consider the following commutative diagram with rows of
distinguished triangles
$$\xymatrix@R=0.3cm@C=0.5cm{\Omega^{-(m-2)}Y_m^{m-1}\ar[r]\ar[d] & \Omega^{-(m-1)}Y_{m+1}^m\ar[r]\ar@{=}[d]&
\Omega^{-(m-1)}Y_{m+1}^{m-1}\ar[r]\ar@{.>}[d]& \Omega^{-(m-1)}Y_m^{m-1}\ar[d]\\
\vdots\ar[d]& \vdots\ar@{=}[d]& \vdots\ar@{.>}[d]& \vdots\ar[d] \\
\Omega^{-(m-2)}Y_m^1 \ar[r]\ar[d]&
\Omega^{-(m-1)}Y_{m+1}^m\ar[r]\ar@{=}[d]&
\Omega^{-(m-1)}Y_{m+1}^{1}\ar[r]\ar@{.>}[d]&\Omega^{-(m-1)}Y_m^1\ar[d]\\
\Omega^{-(m-1)}X_m \ar[r]\ar[d] &
\Omega^{-(m-1)}Y_{m+1}^m\ar[r]\ar@{=}[d]&
\Omega^{-m}X_{m+1}\ar[r]\ar@{.>}[d] & \Omega^{-m}X_m\ar[d]\\
\Omega^{-(m-1)}Y_n^m
\ar[r]\ar[d]&\Omega^{-(m-1)}Y_{m+1}^m\ar[r]\ar@{=}[d] &
\Omega^{-m}Y_{n}^{m+1}\ar[r]\ar@{.>}[d] &\Omega^{-m}Y_n^m\ar[d]\\
\vdots\ar[d]& \vdots\ar@{=}[d]& \vdots\ar@{.>}[d] & \vdots\ar[d] \\
\Omega^{-(m-1)}Y_{m+2}^m\ar[r] & \Omega^{-(m-1)}Y_{m+1}^m\ar[r]
&\Omega^{-m}Y_{m+2}^{m+1}\ar[r] & \Omega^{-m}Y_{m+2}^m}$$ where the
$l$-th row ($1\le l\le m-1$) is from the fourth column of $(3.11)$
by taking $j=m+1, i=m, k=m-l$, and then applying $\Omega^{-(m-1)}$;
the $m$-th row is from the first row of $(3.11)$ by taking $j=m+1,
i=m$, and then applying $\Omega^{-(m-1)}$; the $l$-th row ($m+1\le
l\le n-1$) is from the fourth column of $(3.11)$ by taking
$j=n+m+1-l, i=m+1, k=m$, and then applying $\Omega^{-m}$. By the
definition of the rotation (cf. (3.4)), this proves (3.12) for
$m+1$.

\vskip10pt

\subsection{} For  $X_{(\underline{\phi_i})}\in {\rm Mor}_n(A\mbox{-}\underline {\rm
mod})$, define  $\Omega^{-1} X_{(\underline {\phi_i})}$ to be
$\left(\begin{smallmatrix} \Omega^{-1} X_1\\ \vdots\\
\Omega^{-1} X_n\end{smallmatrix}\right)_{(\Omega^{-1} \underline
{\phi_i})}\in {\rm Mor}_n(A\mbox{-}\underline {\rm mod}).$

\vskip10pt

\begin{lem} \label{rotn+1} \ Let $X_{(\phi_i)}\in {\rm Mor}_n(A)$. Then \ ${\bf Rot}^{j(n+1)} X_{(\underline{\phi_i})}
= \Omega^{-j(n-1)}X_{(\underline{\phi_i})}, \ \forall \ j\ge 1.$
\end{lem} \noindent{\bf Proof.} \  By taking $m=n$ in $(3.12)$, we
get $${\bf Rot}^n
X_{(\underline{\phi_i})}=(\Omega^{-(n-2)}Y_n^{n-1}\rightarrow\Omega^{-(n-2)}Y_n^{n-2}\rightarrow\cdots\rightarrow
\Omega^{-(n-2)}Y_n^{1}\rightarrow\Omega^{-(n-1)}X_n).$$ We have the
following commutative diagram with rows being distinguished
triangles
$$\xymatrix@R=0.5cm@C=1cm{\Omega^{-(n-2)}Y_n^{n-1}\ar[r]\ar[d] & \Omega^{-(n-1)}X_n\ar[rrr]^-{(-1)^{n-1}\Omega^{-(n-1)}\underline{\phi_{n-1}}}\ar@{=}[d] &
& &
\Omega^{-(n-1)}X_{n-1}\ar[r]\ar@{.>}[d]^{\Omega^{-(n-1)}\underline{\phi_{n-2}}}& \Omega^{-(n-1)}Y_n^{n-1}\ar[d]\\
\Omega^{-(n-2)}Y_n^{n-2}\ar[r]\ar[d]&
\Omega^{-(n-1)}X_n\ar[rrr]^-{(-1)^{n-1}\Omega^{-(n-1)}\underline{\phi_{n-2}}\
\underline{\phi_{n-1}}}\ar@{=}[d] & & &
\Omega^{-(n-1)}X_{n-2}\ar[r]\ar@{.>}[d] &\Omega^{-(n-1)}Y_n^{n-2}\ar[d]\\
\vdots\ar[d] & \vdots \ar@{=}[d] & & &\vdots\ar@{.>}[d]^{\Omega^{-(n-1)}\underline{\phi_{2}}} & \vdots \ar[d]\\
\Omega^{-(n-2)}Y_n^{2}\ar[r]\ar[d]& \Omega^{-(n-1)}X_n
\ar[rrr]^-{(-1)^{n-1}\Omega^{-(n-1)}\underline{\phi_2}\cdots\underline{\phi_{n-1}}}\ar@{=}[d]&&&
\Omega^{-(n-1)}X_{2}\ar[r]\ar@{.>}[d]^{\Omega^{-(n-1)}\underline{\phi_{1}}} &\Omega^{-(n-1)}Y_n^{2}\ar[d]\\
\Omega^{-(n-2)}Y_n^{1}\ar[r]&
\Omega^{-(n-1)}X_n\ar[rrr]^-{(-1)^{n-1}\Omega^{-(n-1)}\underline{\phi_1}\cdots\underline{\phi_{n-1}}}
&&& \Omega^{-(n-1)}X_{1}\ar[r]&\Omega^{-(n-1)}Y_n^{1}}$$ where the
$l$-th row ($1\le l\le n-1$) is from the first row of $(3.11)$ by
taking $j=n, i=n-l$, and then applying $\Omega^{-(n-1)}$ (note that
$(-1)^{(n-1)}$ arises from applying $\Omega^{-(n-1)}$). Using (3.4)
we get that ${\bf Rot}^{n+1} X_{(\underline{\phi_i})}$ is
$$\begin{array}{l}
(\xymatrix{\Omega^{-(n-1)}X_n\ar[rrr]^-{(-1)^{n-1}\Omega^{-(n-1)}\underline{\phi_{n-1}}}
& & &
\Omega^{-(n-1)}X_{n-1}\ar[rr]^-{\Omega^{-(n-1)}\underline{\phi_{n-2}}}
&& \cdots\ar[rr]^-{\Omega^{-(n-1)}\underline{\phi_{1}}}
&& \Omega^{-(n-1)}X_{1}}) \\
\cong
(\xymatrix{\Omega^{-(n-1)}X_n\ar[rr]^-{\Omega^{-(n-1)}\underline{\phi_{n-1}}}
& &
\Omega^{-(n-1)}X_{n-1}\ar[rr]^-{\Omega^{-(n-1)}\underline{\phi_{n-2}}}
&& \cdots\ar[rr]^-{\Omega^{-(n-1)}\underline{\phi_{1}}} &&
\Omega^{-(n-1)}X_{1}})\\
= \Omega^{-(n-1)}X_{(\underline{\phi_i})}.\end{array}$$ From this
and induction the assertion follows. $\s$

\vskip 10pt

\subsection{}\ Since $\tau$ is a
triangle functor, by construction we see $${\bf Rot} \ \tau
X_{(\underline {\phi_i})}\cong \tau \ {\bf Rot} X_{(\underline
{\phi_i})}.\eqno(3.13)$$ We have the following important result.

\begin{thm} \label{Maintheorem} \ Let $A$ be a selfinjective algebra, $X_{(\phi_i)}\in\mathcal{S}_n(A)$. Then
there are the following isomorphisms in
${\rm Mor}_n(A\mbox{-}\underline {\rm mod})$

\vskip10pt

$(i)$ \  $\underline{\tau_{\mathcal{S}}^{j}X_{(\phi_i)}}\cong
\tau^{j} \ {\bf Rot}^{j}X_{(\underline{\phi_i})}$ for $j\ge 1$. In
particular, $\underline{\tau_{\mathcal{S}}X_{(\phi_i)}}\cong \tau \
\underline{{\bf Cok}X_{(\phi_i)}}$.

\vskip10pt

$(ii)$ \ $\underline {\tau_{\mathcal{S}}^{s(n+1)}X_{(\phi_i)}}\cong
\tau^{s(n+1)} \ \Omega^{-s(n-1)}X_{( \underline{\phi_i} )}, \
\forall \ s\ge 1$.
\end{thm}

\vskip5pt

\noindent{\bf Proof.} \ $(i)$  \ First, we claim that there are the
following isomorphisms in ${\rm Mor}_n(A\mbox{-}\underline {\rm
mod})$:
$$\underline {({\bf Cok \ Mimo} \ \tau )^j \ Y_{(\psi_i)}} \cong
 \tau^j \ {\bf Rot}^j  Y_{(\underline {\psi_i})}, \ \ \forall \
Y_{(\psi_i)}\in{\rm Mor}_n(A\mbox{-}{\rm mod}), \ \forall \ j\ge
1.\eqno(3.14)$$ In fact, by Lemma 3.1 and induction we have
$$\begin{array}{ll}
\underline {({\bf Cok \ Mimo} \ \tau )^j \ Y_{(\psi_i)}} & \cong
\underline {{\bf Cok \ Mimo} \ \tau \ ({\bf Cok \ Mimo} \ \tau
)^{j-1} \ Y_{(\psi_i)}} \\[1em] &
 \overset{_{Lemma \ 3.1}}\cong
{\bf Rot} \ \underline {\tau \ ({\bf Cok \ Mimo} \ \tau
)^{j-1} \ Y_{(\psi_i)}} \\[1em]
& \overset{(2.3)}\cong {\bf Rot} \ \tau \ \underline {({\bf Cok \
Mimo} \ \tau
)^{j-1} \ Y_{(\psi_i)}} \\[1em]
& \overset{_{Induction}}\cong  {\bf Rot} \ \tau^{j} \ {\bf Rot}^{j-1}   Y_{(\underline {\psi_i})} \\[1em]
&  \overset{(3.13)}\cong \tau^{j} \ {\bf Rot}^{j} Y_{(\underline
{\psi_i})}.
\end{array}$$

\vskip5pt

Now, we have the following isomorphisms in ${\rm
Mor}_n(A\mbox{-}\underline {\rm mod})$:
$$\begin{array}{ll}
\underline{\tau_{\mathcal{S}}^{j}X_{(\phi_i)}}& \overset{_{Theorem \
2.4}}\cong
\underline {({\bf Mimo} \ \tau \ {\bf Cok})^{j}X_{(\phi_i)}}\\[1em]
& \cong \underline {{\bf Mimo} \ \tau \ ({\bf Cok \ Mimo} \ \tau)^{j-1} \ {\bf Cok}X_{(\phi_i)}} \\[1em]
  & \overset{(2.3)}\cong \tau \ \underline {({\bf Cok \ Mimo} \ \tau)^{j-1} \ {\bf Cok}X_{(\phi_i)}}\\[1em]
&  \overset{(3.14)}\cong \tau^{j} \ {\bf Rot}^{j-1} \ \underline {{\bf Cok}X_{(\phi_i)}} \\[1em]
& \cong\tau^{j} \ {\bf Rot}^{j-1} \ \underline {{\bf Cok \ Mimo}X_{(\phi_i)}}\\[1em]
& \overset{_{Lemma \ 3.1}}\cong \tau^{j} \ {\bf
Rot}^{j}X_{(\underline{\phi_i})},
\end{array}$$
where we have used ${\bf Mimo} X_{(\phi_i)} = X_{(\phi_i)}$ since
$X_{(\phi_i)}\in \mathcal{S}_n(A)$.

\vskip5pt

$(ii)$ \ This follows from Lemma 3.3 and $(i)$ by taking $j =
s(n+1)$. $\s$

\subsection{} Now we pass from ${\rm Mor}_n(A\mbox{-}\underline {\rm
mod})$ to ${\rm Mor}_n(A\mbox{-}{\rm mod}).$ Before stating the main
result, we need a notation. Let $X_{(\phi_i)}\in {\rm Mor}_n(A)$.
For positive integers  $r$ and $t$, the object $\tau^r \ \Omega^{-t}
X_{(\underline{\phi_i})}\in {\rm Mor}_n(A\mbox{-}\underline {\rm
mod})$ is already defined (cf. 2.4 and 3.3). As in 2.4, we consider
the full subcategory of ${\rm Mor}_n(A)$ given by
$$\{Y_{(\psi_i)} = \left(\begin{smallmatrix} \tau^r \ \Omega^{-t} X_1\\ \vdots\\
\tau^r \ \Omega^{-t} X_n\end{smallmatrix}\right)_{(\psi_i)}\in {\rm
Mor}_n(A) \ | \ Y_{(\underline {\psi_i})} \cong  \tau^r \
\Omega^{-t} X_{(\underline {\phi_i})}\}.$$ Any object in this
subcategory will be denoted by $\tau^r \ \Omega^{-t} X_{({\phi_i})}$
(we emphasize that this convention will cause no confusions). So, we
have $\underline {\tau^r \ \Omega^{-t} X_{({\phi_i})}} \cong \tau^r
\ \Omega^{-t} X_{(\underline {\phi_i})}.$ By Lemma 1.6$(i)$, ${\bf
Mimo} \ \tau^r \ \Omega^{-t} X_{(\phi_i)}\in \mathcal{S}_n(A)$ is a
well-defined object, and there are the following isomorphisms  in
${\rm Mor}_n(A\mbox{-}\underline {\rm mod})$$$\underline {{\bf Mimo}
\ \tau^r \ \Omega^{-t} X_{(\phi_i)}}\cong \underline {\tau^r \
\Omega^{-t} X_{({\phi_i})}} \cong \tau^r \ \Omega^{-t}
X_{(\underline {\phi_i})}.\eqno(3.15)$$
\begin{thm}\label{Maintheoremcor}
\ Let $A$ be a selfinjective algebra, and  $X_{(\phi_i)}\in
\mathcal{S}_n(A)$. Then we have
$$\tau_{\mathcal{S}}^{s(n+1)}X_{(\phi_i)}\cong {\bf Mimo} \
\tau^{s(n+1)} \ \Omega^{-s(n-1)}X_{(\phi_i)}, \ \ s\ge
1.\eqno(3.16)$$
\end{thm}
\noindent{\bf Proof.} \ By Theorem \ref{Maintheorem}$(ii)$ we have
$$\underline{\tau_{\mathcal{S}}^{s(n+1)} X_{(\phi_i)}}\cong
\tau^{s(n+1)} \ \Omega^{-s(n-1)}X_{(\underline{\phi_i})} \cong
\underline{{\bf Mimo} \
\tau^{s(n+1)}\Omega^{-s(n-1)}X_{(\phi_i)}}.$$ Since ${\bf Mimo} \
\tau^{s(n+1)} \
\Omega^{-s(n-1)}X_{(\phi_i)}\in\mathcal{S}_n(A)_{\mathcal{I}}$ (cf.
Remark $(ii)$ in 1.4), the assertion follows from Corollary
\ref{functorF}. $\s$

\subsection{} We apply Theorem 3.4 to the algebra $\A(m, t)$, which is defined
below. Let $\Bbb Z_m$ be the cyclic quiver with vertices indexed by
the cyclic group $\Bbb Z/m\Bbb Z$ of order $m$, and with arrows
$a_i: \ i \longrightarrow i+1, \ \forall \ i\in \Bbb Z/m\Bbb Z$. Let
$k\Bbb Z_m$ be the path algebra of the quiver $\Bbb Z_m$, $J$ the
ideal generated by all arrows, and $\A(m, t): =k\Bbb Z_m/J^t$ with
$m\ge1, \ t\ge 2$.  Any connected selfinjective Nakayama algebra
over an algebraically closed field is Morita equivalent to $\A(m,
t),$ $m\ge 1, \ t\ge 2$.  Note that $\A(m, t)$ is a Frobenius
algebra of finite representation type, and that $\A(m, t)$ is
symmetric if and only if $m\mid (t-1)$. For the Auslander-Reiten
sequences of $\A(m, t)$ see [ARS], p.197. In the stable category
$\A(m, t)$-\underline {mod}, we have the following information on
the orders of $\tau$ and $\Omega$ (see 5.1 in [CZ])
$$o(\tau) = m; \ \ \ \ \ \ \ \ o(\Omega) =
\begin{cases} m, & \ t = 2;\\ \frac{2m}{(m, t)}, & \ t\ge 3, \end{cases}\eqno(3.17)$$
where $(m, t)$ is the g.c.d  of $m$ and $t$. By (3.16) and (3.17) we
get the following

\begin{cor}\label{Mainresultforcua} \ For an indecomposable
nonprojective object $X_{(\phi_i)}\in\mathcal{S}_n(\A(m,t))$,
$m\ge1, \ t\ge 2$, there are the following isomorphisms:

\vskip5pt

$(i)$ \ If $n$ is odd, then
$\tau_{\mathcal{S}}^{m(n+1)}X_{(\phi_i)}\cong X_{(\phi_i)}$;

\vskip5pt

$(ii)$ \ If $n$ is even, then
$\tau_{\mathcal{S}}^{2m(n+1)}X_{(\phi_i)}\cong X_{(\phi_i)}$.
\end{cor}
\subsection{Example.} Let $A=kQ/\langle\delta\alpha, \beta\gamma,
\alpha\delta-\gamma\beta\rangle$, where $Q$ is the quiver $
\xymatrix@=0.5cm{2\bullet\ar@<1ex>[r]^{\alpha}& 1
\bullet\ar@<1ex>[r]^{\beta}\ar@<1ex>[l]^{\delta}&3
\bullet\ar@<1ex>[l]^{\gamma}}$ Then $A$ is selfinjective with
$\tau\cong\Omega^{-1}$ and $\Omega^6 \cong {\rm id}$ on the object
of $A$-\underline{mod}. The Auslander-Reiten quiver of $A$ is
$$\xymatrix@R=0.2cm@C=0.4cm{& & & {\begin{smallmatrix}3\\1\\2\end{smallmatrix}}\ar[dr] & & &
{\begin{smallmatrix}&1&\\2 & & 3\\& 1&\end{smallmatrix}}\ar[ddr]&\\
3\ar[dr] & &
{\begin{smallmatrix}1\\2\end{smallmatrix}}\ar[dr]\ar[ur]\ar@{.>}[ll]
& &
{\begin{smallmatrix}3\\1\end{smallmatrix}}\ar[dr]\ar@{.>}[ll] & & 2\ar[dr]\ar@{.>}[ll] &\\
& {\begin{smallmatrix}& 1 & \\2& &3 \end{smallmatrix}}\ar[dr]\ar[ur]
& & 1\ar[dr]\ar[ur]\ar@{.>}[ll] & & {\begin{smallmatrix}2& &3\\ & 1&
\end{smallmatrix}}\ar[dr]\ar[ur]\ar[uur]\ar@{.>}[ll]& &
{\begin{smallmatrix}& 1 & \\2& &3 \end{smallmatrix}}\ar@{.>}[ll]\\
2\ar[ur] & &
{\begin{smallmatrix}1\\3\end{smallmatrix}}\ar[dr]\ar[ur]\ar@{.>}[ll]
& &
{\begin{smallmatrix}2\\1\end{smallmatrix}}\ar[ur]\ar@{.>}[ll]& & 3\ar[ur]\ar@{.>}[ll] &\\
& & &  {\begin{smallmatrix}2\\1\\3\end{smallmatrix}}\ar[ur] & & & &
}$$ Let $X_{(\phi_i)}$ be an indecomposable nonprojective object in
$\mathcal{S}_n(A)$. By (3.16), for $s\ge 1$ we have
$\tau_{\mathcal{S}}^{s(n+1)}X_{(\phi_i)}\cong {\bf Mimo} \
\tau^{s(n+1)} \ \Omega^{-s(n-1)}X_{(\phi_i)} \cong {\bf Mimo} \
\Omega^{-2sn}X_{(\phi_i)}$ in $\mathcal{S}_n(A)$. Then by  Remark
$(i)$ in 1.4 we get

\vskip 5pt

$(i)$ \ if $n\equiv 0, \ \mbox{or} \ 3 \ ({\rm mod} 6)$, then
$\tau_{\mathcal{S}}^{n+1}X_{(\phi_i)}\cong X_{(\phi_i)}$; and

\vskip5pt

$(ii)$ \ if $n\equiv \pm1, \ \mbox{or} \ \pm2 \ ({\rm mod} 6)$, then
$\tau_{\mathcal{S}}^{3(n+1)}X_{(\phi_i)}\cong X_{(\phi_i)}$.

\vskip10pt

\section{\bf Serre functors of stable monomorphism categories}

Throughout this section, $A$ is a finite-dimensional selfinjective
algebra over a field. We study the periodicity of the Serre functor
$F_{\mathcal S}$ on the objects of the stable monomorphism category
$\underline{\mathcal{S}_n(A)}$. In particular, $F_{\mathcal
S}^{m(n+1)}X\cong X$ for $X\in \underline{\mathcal{S}_n(\A(m, t))}$.

\vskip 10pt

\subsection {}  Let $\mathcal{A}$ be a Hom-finite Krull-Schmidt triangulated
$k$-category with suspension functor $[1]$. For the Auslander-Reiten
triangles we refer to [H]. In an Auslander-Reiten triangle $X
{\rightarrow} Y {\rightarrow} Z {\rightarrow} X[1]$, the
indecomposable object $X$ is uniquely determined by $Z$. Write $X =
\widetilde{\tau_\mathcal A} Z$, and extend the action of
$\widetilde{\tau_\mathcal A}$ to arbitrary objects, and put
$\widetilde{\tau_{\mathcal{A}}}0 = 0$. In general,
$\widetilde{\tau_\mathcal A}$ is {\em not} a functor. By Theorem
I.2.4 of [RV], $\mathcal{A}$ has a Serre functor $F$ if and if
$\mathcal{A}$ has Auslander-Reiten triangles; if this is the case,
$F$ and $[1]\widetilde{\tau_\mathcal A}$ coincide on the objects of
$\mathcal A$. If $X\stackrel{f} {\rightarrow} Y \stackrel{g}
{\rightarrow} Z \stackrel{h} {\rightarrow} X[1]$ is an
Auslander-Reiten triangle, then so is $X[1]\stackrel{-f[1]}
{\longrightarrow} Y[1] \stackrel{-g[1]} {\longrightarrow} Z[1]
\stackrel{-h[1]} {\longrightarrow} X[2]$. It follows that $$FZ\cong
([1]\widetilde{\tau_\mathcal A}) Z \cong (\widetilde {\tau_\mathcal
A}[1])Z, \ \forall \ Z\in\mathcal A. \eqno(4.1)$$

\subsection{} Since $A$ is self-injective, by Corollary 4.1$(ii)$ of \cite{Z},
$\mathcal{S}_n(A)$ is exactly the category of Gorenstein projective
$T_n(A)$-modules, hence it is a Frobenius category whose
projective-injective objects are exactly all the projective
$T_n(A)$-modules. Thus, the stable category
$\underline{\mathcal{S}_n(A)}$ of $\mathcal{S}_n(A)$ modulo
projective objects is a Hom-finite Krull-Schmidt triangulated
category with suspension functor $\Omega^{-1}_{\mathcal{S}} =
\Omega^{-1}_{\mathcal{S}_n(A)}$. Since $\mathcal{S}_n(A)$ has
Auslander-Reiten sequences, it follows that
$\underline{\mathcal{S}_n(A)}$ has Auslander-Reiten triangles, and
hence it has a Serre functor $F_{\mathcal{S}} =
F_{\mathcal{S}_n(A)}$, which coincides with
$\Omega_{\mathcal{S}}^{-1}\widetilde{\tau_{\mathcal{S}}}$ on the
objects of $\underline{\mathcal{S}_n(A)}$.

In order to make the following computation more clear, we denote by
$Q: \mathcal S_n(A) \rightarrow \underline{\mathcal{S}_n(A)}$ the
natural functor. Then
$$\widetilde{\tau_{\mathcal{S}}} \ Q Z =  Q \ \tau_{\mathcal{S}}Z =
\tau_{\mathcal{S}}Z, \ \forall Z\in\mathcal S_n(A).\eqno(4.2)$$

\subsection{} \ In 2.6 we have considered the canonical functor $W:
\mathcal{S}_n(A)_{\mathcal{I}}\rightarrow {\rm
Mor}_n(A\mbox{-}\underline {\rm mod})$ given by $X_{(\phi_i)}\mapsto
X_{(\underline {\phi_i})}$.

\vskip10pt

\begin{lem} \ The functor $W$ induces a functor
$\widetilde{W}: \underline{\mathcal S_n(A)}\rightarrow {\rm
Mor}_n(A\mbox{-} \underline{\rm mod})$ satisfying $\widetilde{W} \
Q|_{\mathcal{S}_n(A)_{\mathcal{I}}} = W$; and $\widetilde{W}$
reflects isomorphisms.
\end{lem}
\noindent{\bf Proof.} The definition of $\widetilde{W}$ is clear by
the requirement $\widetilde{W}  Q|_{\mathcal{S}_n(A)_{\mathcal{I}}}
= W$. We need to check that it is well-defined. If $X_{(\phi_i)}$
and $Y_{(\psi_i)}$ are indecomposable and nonprojective in $\mathcal
S_n(A)$, and $X_{(\phi_i)}\cong Y_{(\psi_i)}$ in
$\underline{\mathcal S_n(A)}$, then $X_{(\phi_i)}\cong Y_{(\psi_i)}$
in ${\mathcal S_n(A)}$, and hence $X_{(\underline {\phi_i})}\cong
Y_{(\underline {\psi_i})}$ in ${\rm Mor}_n(A\mbox{-} \underline{\rm
mod})$, i.e., $\widetilde{W}$ is well-defined on objects. For any
morphism $f = {(f_i)}: X_{({\phi_i})}\rightarrow Y_{({\psi_i})}$ in
${\mathcal S_n(A)}$ which factors through a projective object of
$\mathcal S_n(A)$, by Lemma 1.1$(iii)$, the morphism $(\underline
{f_i}) =0$ in ${\rm Mor}_n(A\mbox{-}\underline {\rm mod})$. Thus
$\widetilde{W}$  is well-defined.

Assume that $\widetilde{W} X\cong \widetilde{W} Y$ in ${\rm
Mor}_n(A\mbox{-}\underline{ {\rm mod}})$ for $X, Y\in
\underline{{\mathcal S_n(A)}}.$ We may write $X = QX', \ Y = QY'$
with $X', Y' \in {\mathcal S_n(A)}_{\mathcal I}.$ Then $WX'\cong
WY'$ in ${\rm Mor}_n(A\mbox{-}\underline{ {\rm mod}})$. By Corollary
2.6, $W$ reflects isomorphisms, thus $X'\cong Y'$ in ${\mathcal
S_n(A)}_{\mathcal I}$, and hence $X\cong Y$ in $\underline{\mathcal
S_n(A)}$. $\s$

\vskip10pt

\begin{lem} \ For $X = X_{(\phi_i)}\in \underline{\mathcal{S}_n(A)}$, we
have the following isomorphism in ${\rm Mor}_n(A\mbox{-}\underline
{\rm mod})$$$\widetilde{W} \ \Omega_{\mathcal S}X\cong  \Omega \
\widetilde{W} X.\eqno(4.3)$$
\end{lem}
\noindent{\bf Proof.}  Let $0\rightarrow \Omega_{\mathcal S}X
\rightarrow P \rightarrow X \rightarrow 0$ be an exact sequence in
${\mathcal{S}_n(A)}$ with $P$ projective. Taking the $i$-th branches
we see that $(\Omega_{\mathcal S}X)_i =  \Omega X_i \oplus P'_i$ for
some projective $A$-module $P_i'$. It follows that
$(\Omega_{\mathcal S}X)_i =  \Omega X_i$ in $A$-\underline{mod}.
Write $\Omega_{\mathcal S}X$ as $(\Omega_{\mathcal S}X)_{(\psi_i)}$.
By the following commutative diagram with exact columns
$$\xymatrix@R=0.4cm@C=0.6cm{0\ar[d] & 0\ar[d] &  & 0\ar[d] & 0\ar[d]\\
(\Omega_{\mathcal S}X)_n\ar[r]^-{\psi_{n-1}}\ar[d] &
(\Omega_{\mathcal S}X)_{n-1}\ar[r]\ar[d] & \cdots\ar[r] &
(\Omega_{\mathcal S}X)_2\ar[r]^-{\psi_{1}}\ar[d] & (\Omega_{\mathcal S}X)_1\ar[d] \\
P_n\ar[r]\ar[d] & P_{n-1}\ar[r]\ar[d] & \cdots\ar[r] & P_2\ar[r]\ar[d] & P_1\ar[d]\\
X_n\ar[r]^-{\phi_{n-1}}\ar[d] & X_{n-1}\ar[r]\ar[d] & \cdots\ar[r] & X_2\ar[r]^-{\phi_{1}}\ar[d] & X_1\ar[d]\\
0 & 0 &  & 0 & 0 }$$ we see $\underline {\psi_i} = {\Omega
{\underline {\phi_i}}}$, and hence the assertion follows.  $\s$

\subsection{} For positive integers $a$ and $b$,  let $[a, b]$ denote  the l.c.m
of $a$ and $b$. The main result of this section is

\vskip10pt

\begin{thm} \ Let $A$ be a selfinjective algebra, and $F_{\mathcal{S}}$  be the Serre functor of
$\underline{\mathcal{S}_n(A)}$. Then we have an isomorphism in
$\underline{\mathcal{S}_n(A)}$ for $X_{({\phi_i})}\in
\underline{\mathcal{S}_n(A)}$ and for $s\ge 1$
$$F_{\mathcal{S}}^{s(n+1)} X_{(\phi_i)}\cong {\bf Mimo} \ \tau^{s(n+1)} \
\Omega^{-2sn} X_{({\phi_i})}.\eqno(4.4)$$

Moreover, if $d_1$ and $d_2$ are positive integers such that
$\tau^{d_1}M \cong M$ and $\Omega^{d_2}M\cong M$ for each
indecomposable nonprojective $A$-module $M,$ then
$F^{N(n+1)}_{\mathcal{S}}X_{(\phi_i)} \cong X_{(\phi_i)}$, where $N
= [\frac{d_1}{(n+1, d_1)}, \frac{d_2}{(2n, d_2)}]$.
\end{thm}

\noindent{\bf Proof.} \ We have isomorphisms in ${\rm
Mor}_n(A\mbox{-} \underline{\rm mod})$ for $s\ge 1$:
$$\begin{array}{ll} \widetilde{W} \ F_{\mathcal{S}}^{s(n+1)} X_{(\phi_i)}
& \overset{(4.1)}{\cong} \widetilde{W} \ \Omega_{\mathcal
  S}^{-s(n+1)} \ \widetilde{\tau_{\mathcal S}}^{s(n+1)} X_{(\phi_i)}\\[0.6em]
&\overset{(4.3)}{\cong}  \Omega^{-s(n+1)} \ \widetilde{W} \ \widetilde{\tau_{\mathcal S}}^{s(n+1)} X_{(\phi_i)}\\[0.6em]
&\overset{(4.2)}{\cong}  \Omega^{-s(n+1)} \ \widetilde{W} \ Q \ \tau_{\mathcal S}^{s(n+1)} X_{(\phi_i)}\\[0.6em]
& \overset{(3.16)}{\cong}    \Omega^{-s(n+1)} \ \widetilde{W} \ Q \ {\bf Mimo} \ \tau^{s(n+1)} \ \Omega^{-s(n-1)}X_{(\phi_i)}\\[0.6em]
& \overset{_{Lemma \ 4.1}}{\cong} \Omega^{-s(n+1)} \ W \ {\bf Mimo} \ \tau^{s(n+1)} \ \Omega^{-s(n-1)}X_{(\phi_i)}\\[0.6em]
& \cong \Omega^{-s(n+1)} \
\underline{\tau^{s(n+1)} \ \Omega^{-s(n-1)} X_{(\phi_i)}}\\[0.6em]
& \overset{(3.15)}\cong \Omega^{-s(n+1)} \ \tau^{s(n+1)} \ \Omega^{-s(n-1)}X_{(\underline{\phi_i})}\\[0.6em]
  & \cong  \tau^{s(n+1)} \ \Omega^{-2sn} X_{(\underline{\phi_i})} \\[0.6em]
  & \overset{(3.15)}\cong  \widetilde{W} \ {\bf Mimo} \ \tau^{s(n+1)} \ \Omega^{-2sn}
  X_{({\phi_i})}.
  \end{array}$$
Now (4.4) follows from Lemma 4.1. Since $d_1|N(n+1), \ d_2|(2Nn)$,
taking $s = N$ in (4.4) we get \
$F^{N(n+1)}_{\mathcal{S}}X_{(\phi_i)} \cong {\bf Mimo} \
\tau^{N(n+1)} \ \Omega^{-2Nn} X_{({\phi_i})} \cong {\bf Mimo}
X_{({\phi_i})} = X_{({\phi_i})}. $ $\s$ \vskip10pt

Note that the conditions on $\tau$ and $\Omega$ in Theorem 4.3 hold
in particular for representation-finite selfinjective algebras.

 \vskip10pt
\subsection{} Applying Theorem 4.3 to the selfinjective Nakayama algebras $\A(m, t)$,
we get

\begin{cor} \ Let $F_{\mathcal{S}}$  be the Serre functor of
$\underline{\mathcal{S}_n(\A(m, t))}$ with $m\ge 1, \ t\ge 2$, and
$X$ be \\[0.8em] an arbitrary object in $\underline{\mathcal{S}_n(\A(m,
t))}$. Then

\vskip5pt

$(i)$ \ If $t=2$, then $F^{N(n+1)}_{\mathcal{S}}X \cong X$, where $N
= \frac{m}{(m, n-1)}.$

\vskip5pt

$(ii)$ \  If $t\ge 3$, then $F^{N(n+1)}_{\mathcal{S}}X \cong X$,
where $N = \frac{m}{(m, t, n+1)}$.
\end{cor}

\noindent{\bf Proof.} \ $(i)$ \ In this case $\tau = \Omega$,  and
$o(\tau) = m = o(\Omega)$, by (3.17). Put $N = \frac{m}{(m, n-1)}.$
It follows from (4.4) that
\begin{align*}F_{\mathcal{S}}^{N(n+1)} X& \cong {\bf Mimo} \
\tau^{N(n+1)} \ \Omega^{-2Nn} X\cong {\bf Mimo} \ \Omega^{-N(n-1)} X
\\ &\cong {\bf Mimo} \ \Omega^{-m \frac {n-1}{(m, n-1)}} X \cong {\bf Mimo} X =
X.\end{align*}

\vskip5pt

$(ii)$ \ This follows from Theorem 4.3 by taking $d_1 = m$ and $d_2
= \frac{2m}{(m, t)}$. By a computation in elementary number theory,
we get $$ \ \ \ \ \ \ \ \ \ \ [\frac{d_1}{(n+1, d_1)},
\frac{d_2}{(2n, d_2)}]=[\frac{m}{(n+1, m)}, \frac{\frac{2m}{(m,
t)}}{(2n, \frac{2m}{(m, t)})}]=\frac{m}{(m, t, n+1)}. \ \ \ \ \  \ \
\ \ \ \ \ \s$$

\vskip10pt

\section{\bf Appendix 1: Proofs of Lemmas 1.5 and 1.6}

We give proofs of Lemmas 1.5 and 1.6.

\vskip10pt

\begin{lem}{\label{isoobjects2}}
\ Let $X_{(\phi_i)}\in{\rm Mor}_n(A)$,  \ $I_2, \cdots, I_n$ be
injective $A$-modules such that $X'_{(\phi'_i)}=\left(\begin{smallmatrix}X_1\oplus I_2\oplus\cdots\oplus I_n\\
\vdots\\X_{n-1}\oplus
I_n\\X_n\end{smallmatrix}\right)_{(\phi'_i)}\in \mathcal{S}_n(A), \
\  \mbox{where}\ \  \phi'_i= \left(\begin{smallmatrix}
\phi_i& 0 & 0  &\cdots &  0 \\
*& 0 & 0 & \cdots & 0 \\
* & 1& 0 & \cdots &  0 \\
 \vdots & \vdots & \vdots &
\ddots & \vdots \\
* & 0& 0 & \cdots & 1
\end{smallmatrix}\right)_{(n-i+1)\times (n-i)}.$
Then $X'_{(\phi'_i)}\cong{\bf Mimo}X_{(\phi_i)}\oplus J$, where $J$
is an injective object of $\mathcal{S}_n(A)$. Moreover,
$J=\bigoplus\limits_{i=1}^{n-1}{\bf m}_i(Q_i)$, where $Q_i$ is an
injective $A$-module such that $Q_i\oplus {\rm IKer}\phi_i\cong
I_{i+1}$, $1\le i\le n-1$.
\end{lem}
\noindent{\bf Proof.} \ It is clear that the morphism $\left(\begin{smallmatrix}(1, 0, \cdots, 0)\\
\vdots\\(1, 0)\\1\end{smallmatrix}\right): X'_{(\phi'_i)}
\twoheadrightarrow X_{(\phi_i)}$ is a right approximation of
$X_{(\phi_i)}$ in $\mathcal{S}_n(A)$ (this can be proved as Lemma
1.3$(i)$, see [Z], Lemma 2.3). By Lemma 1.3$(i)$, there is an object
$J\in \mathcal S_n(A)$ such that $X'_{(\phi'_i)}\cong{\bf
Mimo}X_{(\phi_i)}\oplus J$. Comparing the branches we get $J_n = 0$
and
$$I_{i+1}\oplus \cdots\oplus I_n \cong {\rm IKer}\phi_i \oplus \cdots \oplus {\rm IKer}\phi_{n-1} \oplus J_i,
\ \forall \ 1\le i\le n-1.\eqno(5.1)$$ Put $Q_{n-1}=J_{n-1}$. From
$I_n \cong  {\rm IKer}\phi_{n-1} \oplus J_{n-1}$ we see that
$Q_{n-1}$ is an injective $A$-module. Since $J\in \mathcal S_n(A)$,
$Q_{n-1}$ is a submodule of $J_{n-2}$, thus $J_{n-2} = Q_{n-2}\oplus
Q_{n-1}$. By $I_{n-1}\oplus I_n  \cong  {\rm IKer}\phi_{n-2} \oplus
{\rm IKer}\phi_{n-1}\oplus J_{n-2}$ in (5.1), we see $I_{n-1} \cong
{\rm IKer}\phi_{n-2}\oplus Q_{n-2}$. Repeating this process we see
that $J_i$ is of the form $J_i = Q_{i}\oplus \cdots \oplus Q_{n-1}$
with $Q_i$ being injective $A$-modules, and $Q_i\oplus {\rm
IKer}\phi_i\cong I_{i+1}$, \ $1\le i\le n-1$. Thus $J =
\bigoplus\limits_{i=1}^{n-1}{\bf m}_i(Q_i)$ is an injective object
of $\mathcal{S}_n(A)$. $\s$

\vskip 10pt

Now we can prove Lemma 1.5 (cf. Claim 2 in $\S$4 of [RS2] for the
case of $n=2$).

\vskip 10pt

 {\bf Proof of Lemma 1.5.} \ We just prove $(i)$. Since the
$A$-map $\phi'_i: X_{i+1}\oplus I_{i+2}\oplus\cdots\oplus
I_n\rightarrow X_i\oplus I_{i+1}\oplus I_{i+2}\oplus\cdots\oplus
I_n$ is monic, the restriction of $\phi'_i$ on
$I_{i+2}\oplus\cdots\oplus I_n$ is also monic, and hence it is a
split monomorphism. Hence $X'_{(\phi'_i)}$ is isomorphic to
$$X''_{(\phi''_i)}=\left(\begin{smallmatrix}X_1\oplus I_2\oplus\cdots\oplus I_n\\
\vdots\\X_{n-1}\oplus
I_n\\X_n\end{smallmatrix}\right)_{(\phi''_i)}\in \mathcal{S}_n(A), \
\  \mbox{where}\ \  \phi''_i= \left(\begin{smallmatrix}
\phi_i& 0 & 0  &\cdots &  0 \\
*& 0 & 0 & \cdots & 0 \\
* & 1& 0 & \cdots &  0 \\
 \vdots & \vdots & \vdots &
\ddots & \vdots \\
* & 0& 0 & \cdots & 1
\end{smallmatrix}\right)_{(n-i+1)\times (n-i)}.$$ Then the assertion follows from Lemma
\ref{isoobjects2}. $\s$

\vskip10pt
\begin{lem}{\label{isoobjects1}}
\ Let $X_{(f_i)}, \ X_{(g_i)}\in {\rm Mor}_n(A)$ such that $f_i-g_i$
factors through an injective $A$-module, $1\leq i\leq n-1$, and
$h_i: X_i\rightarrow I_i$ be an injective envelope, $2\leq i\leq n$.
Set $$X'_{(f'_i)}=\left(\begin{smallmatrix}
X_1\oplus I_2\oplus\cdots\oplus I_n\\
\vdots
\\X_{n-1}\oplus I_n \\X_n
\end{smallmatrix}\right)_{(f'_i)}  {\rm where}\ \ \
f'_i = \left(\begin{smallmatrix}
f_i& 0 & 0  &\cdots &  0 \\
h_{i+1}& 0 & 0 & \cdots &  0 \\
0 & 1& 0 & \cdots &  0 \\
 \vdots & \vdots & \vdots &
\ddots & \vdots \\
0 & 0& 0 & \cdots & 1
\end{smallmatrix}\right)_{(n-i+1)\times (n-i)},$$
and
$$X'_{(g'_i)}=\left(\begin{smallmatrix}
X_1\oplus I_2\oplus\cdots\oplus I_n\\
\vdots
\\X_{n-1}\oplus I_n\\X_n
\end{smallmatrix}\right)_{(g'_i)} {\rm where}\ \ \
g'_i = \left(\begin{smallmatrix}
g_i& 0 & 0  &\cdots &  0 \\
h_{i+1}& 0 & 0 & \cdots &  0 \\
0 & 1& 0 & \cdots &  0 \\
 \vdots & \vdots & \vdots &
\ddots & \vdots \\
0 & 0& 0 & \cdots & 1
\end{smallmatrix}\right)_{(n-i+1)\times (n-i)}.$$
Then $X'_{(f'_i)}\cong X'_{(g'_i)}$ in $\mathcal{S}_n(A)$.
\end{lem}
\noindent{\bf Proof.} For $1\leq i\leq n-1$, it is clear that
$f_i-g_i: X_{i+1}\rightarrow X_i$ factors through the injective
envelope $h_{i+1}: X_{i+1}\rightarrow I_{i+1}$, and hence there is
an $A$-map $u_i: I_{i+1}\rightarrow X_i$ such that
$g_i-f_i=u_ih_{i+1}$. The following commutative diagram shows
$X'_{(f'_i)}\cong X'_{(g'_i)}$.
$$\xymatrix@R=0.7cm@C=0.5cm{{X_n}\ar[r]^-{f'_{n-1}}\ar@{=}[d]_{\alpha_n=1_{X_n}} &
{X_{n-1}\oplus I_n}\ar[r]\ar[d]^{\alpha_{n-1}} & \cdots\ar[r] &
{X_2\oplus I_3\oplus\cdots\oplus I_n}\ar[d]^{\alpha_2}\ar[r]^-{f'_1}
& {X_1\oplus I_2\oplus\cdots\oplus I_n}\ar[d]^{\alpha_1}
\\ {X_n}\ar[r]^-{g'_{n-1}} & {X_{n-1}\oplus I_n}\ar[r] &
\cdots\ar[r] & {X_2\oplus I_3\oplus\cdots\oplus I_n}\ar[r]^-{g'_1} &
{X_1\oplus I_2\oplus\cdots\oplus I_n}}$$ where $\alpha_i =
\left(\begin{smallmatrix}
1& u_i & g_iu_{i+1}  &\cdots &  (g_ig_{i+1}\cdots g_{n-2}u_{n-1}) \\
0& 1 & h_{i+1}u_{i+1} & \cdots &  (h_{i+1}g_{i+1}\cdots g_{n-2}u_{n-1}) \\
0 & 0& 1 & \cdots &  (h_{i+2}g_{i+2}\cdots g_{n-2}u_{n-1}) \\
 \vdots & \vdots & \vdots &
\ddots & \vdots \\
0 & 0& 0 & \cdots & h_{n-1}u_{n-1}\\
0 & 0& 0 & \cdots & 1
\end{smallmatrix}\right)_{(n-i+1)\times (n-i+1)}, \ 1\leq i\leq
n$. $\s $

\vskip10pt

\begin{lem}\label{isoobjects}
\ Let $X_{(f_i)}, \ X_{(g_i)}\in {\rm Mor}_n(A)$ such that $f_i-g_i$
factors through an injective $A$-module, $1\leq i\leq n-1$. If each
$X_i$ has no nonzero injective direct summands, then ${\bf
Mimo}X_{(f_i)}\cong{\bf Mimo}X_{(g_i)}$ in $\mathcal{S}_n(A)$.
\end{lem}
\noindent{\bf Proof.} \ Consider $X'_{(f'_i)}$ and $X'_{(g'_i)}$
defined in Lemma \ref{isoobjects1}, which are isomorphic in
$\mathcal{S}_n(A)$. By Lemma \ref{isoobjects2}, there exist
injective $A$-modules $Q_{f, i}$ and $Q_{g, i}$ such that ${\rm
IKer}f_i\oplus Q_{f, i}\cong I_{i+1}\cong{\rm IKer}g_i\oplus
Q_{g,i}, \ 1\leq i\leq n-1,$ and
$${\bf Mimo}X_{(f_i)} \oplus\bigoplus\limits_{i=1}^{n-1}{\bf
m}_i(Q_{f, i})\cong X'_{(f'_i)}\cong X'_{(g'_i)}\cong{\bf
Mimo}X_{(g_i)}\oplus\bigoplus\limits_{i=1}^{n-1}{\bf m}_i(Q_{g,
i}).$$ By Claim 3 in $\S4$ of [RS2], we have ${\rm IKer}f_i\cong{\rm
IKer}g_i, \ 1\le i\le n-1$. Thus $Q_{f,i}\cong Q_{g, i}$, and
$\bigoplus\limits_{i=1}^{n-1}{\bf m}_i(Q_{f,
i})\cong\bigoplus\limits_{i=1}^{n-1}{\bf m}_i(Q_{g, i})$, from which
the assertion follows. $\s$

\vskip10pt

Now we can prove Lemma \ref{MimoIso} (cf. Theorem 4.2 of [RS2] for
the case of $n=2$).

\vskip10pt

{\bf Proof of Lemma \ref{MimoIso}.}\ We only prove $(i)$. If ${\bf
Mimo}X_{(\phi_i)} \cong {\bf Mimo} Y_{(\psi_i)}$ in
$\mathcal{S}_n(A)$, then by the construction of ${\bf Mimo}$,
$X_{(\overline{\phi_i})}\cong Y_{(\overline{\psi_i})}$ in ${\rm
Mor}_n(A\mbox{-}\overline {\rm mod})$. Conversely, assume that
$X_{(\overline{\phi_i})}\cong Y_{(\overline{\psi_i})}$ in ${\rm
Mor}_n(A\mbox{-}\overline {\rm mod})$. Since all $X_i$ and $Y_i$
have no nonzero injective direct summands, there are
$A$-isomorphisms $x_i: X_i\rightarrow Y_i$, such that each
$x_i\phi_i-\psi_ix_{i+1}$ factors through an injective $A$-module,
$1\leq i\leq n-1$. Then each $\phi_i-x_i^{-1}\psi_ix_{i+1}$ factors
through an injective $A$-module, $1\leq i\leq n-1$. By Lemma
\ref{isoobjects}, ${\bf Mimo}X_{(x_i^{-1}\psi_ix_{i+1})}\cong{\bf
Mimo}X_{(\phi_i)}$. Since $Y_{(\psi_i)}\cong
X_{(x_i^{-1}\psi_ix_{i+1})}$, we get ${\bf Mimo}Y_{(\psi_i)}\cong
{\bf Mimo}X_{(\phi_i)}. $ $\s$

\section{\bf Appendix 2: Auslander-Reiten quivers  of some monomorphism categories}

We include the Auslander-Reiten quivers of some
representation-finite monomorphism categories.
\subsection{} By Simson \cite{S}, $\mathcal{S}_{n, 2}$, $\mathcal{S}_{2, 3}$,
$\mathcal{S}_{2, 4}$, $\mathcal{S}_{2, 5}$, $\mathcal{S}_{3, 3}$ and
$\mathcal{S}_{4, 3}$ are the only representation-finite cases among
all $\mathcal{S}_{n, t} = \mathcal{S}_n(k[x]/\langle x^t\rangle), \
n\ge 2, \ t\ge 2$. In \cite{RS3}, Ringel and Schmidmeier give the
Auslander-Reiten quivers of $\mathcal{S}_{2, t}$,  $t=2, 3, 4, 5$.
In the following we give the remaining cases, namely,
$\mathcal{S}_{n, 2}$ \ $(n\ge 3)$, \ $\mathcal{S}_{3, 3}$, and
$\mathcal{S}_{4, 3}$.

$(i)$ \ There are $n$ indecomposable projective objects and
$\frac{n(n+1)}{2}$ indecomposable nonprojective objects in
$\mathcal{S}_{n, 2}$. For the Auslander-Reiten quivers of
$\mathcal{S}_{3, 2}$ see 2.5. The Auslander-Reiten quiver of
$\mathcal{S}_{4, 2}$ is as follows, where $A=k[x]/\langle
x^2\rangle$ and $S$ is the simple $A$-module.
$$\xymatrix@R=0.01cm@C=0.3cm{& {\begin{smallmatrix}A\\0\\0\\0\end{smallmatrix}}\ar[dr] & &
{\begin{smallmatrix}A\\A\\0\\0\end{smallmatrix}}\ar[dr] & &
{\begin{smallmatrix}A\\A\\A\\0\end{smallmatrix}}\ar[dr] & &
{\begin{smallmatrix}A\\A\\A\\A\end{smallmatrix}}\ar[dr] &\\
{\begin{smallmatrix}S\\0\\0\\0\end{smallmatrix}}\ar[ur]\ar[dr] & &
{\begin{smallmatrix}A\\S\\0\\0\end{smallmatrix}}\ar[ur]\ar[dr]\ar@{.>}[ll]
& &
{\begin{smallmatrix}A\\A\\S\\0\end{smallmatrix}}\ar[ur]\ar[dr]\ar@{.>}[ll]
& &
{\begin{smallmatrix}A\\A\\A\\S\end{smallmatrix}}\ar[ur]\ar[dr]\ar@{.>}[ll]
& &
{\begin{smallmatrix}S\\S\\S\\S\end{smallmatrix}}\ar@{.>}[ll]\\
& {\begin{smallmatrix}S\\S\\0\\0\end{smallmatrix}}\ar[ur]\ar[dr] & &
{\begin{smallmatrix}A\\S\\S\\0\end{smallmatrix}}\ar[ur]\ar[dr]\ar@{.>}[ll]
& &
{\begin{smallmatrix}A\\A\\S\\S\end{smallmatrix}}\ar[ur]\ar[dr]\ar@{.>}[ll]
& &
{\begin{smallmatrix}S\\S\\S\\0\end{smallmatrix}}\ar[ur]\ar@{.>}[ll] &\\
& & {\begin{smallmatrix}S\\S\\S\\0\end{smallmatrix}}\ar[ur]\ar[dr] &
&
{\begin{smallmatrix}A\\S\\S\\S\end{smallmatrix}}\ar[ur]\ar[dr]\ar@{.>}[ll]
& &
{\begin{smallmatrix}S\\S\\0\\0\end{smallmatrix}}\ar[ur]\ar@{.>}[ll] & &\\
& & & {\begin{smallmatrix}S\\S\\S\\S\end{smallmatrix}}\ar[ur] & &
{\begin{smallmatrix}S\\0\\0\\0\end{smallmatrix}}\ar[ur]\ar@{.>}[ll]
& & & }$$

$(ii)$ \ Let $A=k[x]/\langle x^3\rangle$. Denote by $M$ and $S$ the
two indecomposable nonprojective $A$-modules, where $S$ is simple.
There are 3 indecomposable projective objects and 24 indecomposable
nonprojective objects in $\mathcal{S}_{3, 3}$, whose
Auslander-Reiten quiver is as follows.
$$\xymatrix@R=0.05cm@C=0.15cm{ & {\begin{smallmatrix}A\\0\\0\end{smallmatrix}}\ar[dr]& & & &
{\begin{smallmatrix}A\\A\\A\end{smallmatrix}}\ar[dr] & & & & \\
{\begin{smallmatrix}M\\0\\0\end{smallmatrix}}\ar[dr]\ar[ur] & &
{\begin{smallmatrix}A\\S\\0\end{smallmatrix}}\ar[dr]\ar@{.>}[ll] & &
{\begin{smallmatrix}A\\A\\M\end{smallmatrix}}\ar[dr]\ar[ur]\ar@{.>}[ll]
& & {\begin{smallmatrix}M\\M\\M\end{smallmatrix}}\ar[dr]\ar@{.>}[ll]
& &
{\begin{smallmatrix}M\\0\\0\end{smallmatrix}}\ar[dr]\ar@{.>}[ll] & \\
& {\begin{smallmatrix}M\\S\\0\end{smallmatrix}}\ar[dr]\ar[ur] & &
{\begin{smallmatrix}A\oplus A \\S\oplus
A\\M\end{smallmatrix}}\ar[dr]\ar[ur]\ar@{.>}[ll] & &
{\begin{smallmatrix}M\\M\\S\end{smallmatrix}}\ar[dr]\ar[ur]\ar@{.>}[ll]
& & {\begin{smallmatrix}S\oplus
A\\M\\M\end{smallmatrix}}\ar[dr]\ar[ur]\ar@{.>}[ll] & &
{\begin{smallmatrix}M\\S\\0\end{smallmatrix}}\ar@{.>}[ll] \\
{\begin{smallmatrix}S\oplus A\\S\oplus
M\\M\end{smallmatrix}}\ar[dr]\ar[ur]\ar[ddr] &&
{\begin{smallmatrix}M\oplus A\\S\oplus
A\\M\end{smallmatrix}}\ar[dr]\ar[ur]\ar[ddr]\ar@{.>}[ll] & &
{\begin{smallmatrix}M\oplus A\\S\oplus
M\\S\end{smallmatrix}}\ar[dr]\ar[ur]\ar[ddr]\ar@{.>}[ll] & &
{\begin{smallmatrix}S\oplus
A\\M\\S\end{smallmatrix}}\ar[dr]\ar[ur]\ar[ddr]\ar@{.>}[ll] & &
 {\begin{smallmatrix}S\oplus A\\S\oplus
M\\M\end{smallmatrix}}\ar[dr]\ar[ur]\ar[ddr]\ar@{.>}[ll]& \\
&{\begin{smallmatrix}A\\M\\M\end{smallmatrix}}\ar[ur] & &
{\begin{smallmatrix}M\\M\\0\end{smallmatrix}}\ar[ur]\ar@{.>}[ll] & &
{\begin{smallmatrix}A\\S\\S\end{smallmatrix}}\ar[ur]\ar@{.>}[ll] & &
{\begin{smallmatrix}S\\S\\0\end{smallmatrix}}\ar[ur]\ar@{.>}[ll] & &
{\begin{smallmatrix}A\\M\\M\end{smallmatrix}}\ar@{.>}[ll]  \\
& {\begin{smallmatrix}S\oplus A\\S\oplus
A\\M\end{smallmatrix}}\ar[uur]\ar[dr]& &
{\begin{smallmatrix}M\\S\\S\end{smallmatrix}}\ar[uur]\ar[dr]\ar@{.>}[ll]
& & {\begin{smallmatrix}S\oplus
A\\M\\0\end{smallmatrix}}\ar[uur]\ar[dr]\ar@{.>}[ll] & &
{\begin{smallmatrix}A\\M\\S\end{smallmatrix}}\ar[uur]\ar[dr]\ar@{.>}[ll]
& & {\begin{smallmatrix}S\oplus A\\S\oplus
A\\M\end{smallmatrix}}\ar@{.>}[ll] \\
{\begin{smallmatrix}A\\A\\S\end{smallmatrix}}\ar[ur]&
&{\begin{smallmatrix}S\\S\\S\end{smallmatrix}}\ar[ur]\ar@{.>}[ll] &
& {\begin{smallmatrix}S\\0\\0\end{smallmatrix}}\ar[ur]\ar@{.>}[ll] &
&{\begin{smallmatrix}A\\M\\0\end{smallmatrix}}\ar[ur]\ar[dr]\ar@{.>}[ll]
& & {\begin{smallmatrix}A\\A\\S\end{smallmatrix}}\ar[ur]\ar@{.>}[ll] &\\
& & & &  & & & {\begin{smallmatrix}A\\A\\0\end{smallmatrix}}\ar[ur]&
& }$$

$(iii)$ \ Let $A, M, S$ be as in $(ii)$. There are 4 indecomposable
projective objects and 80 indecomposable nonprojective objects in
$\mathcal{S}_{4, 3}$, whose Auslander-Reiten quiver looks like.

\vskip15pt

\begin{center}
\includegraphics[scale=0.4]{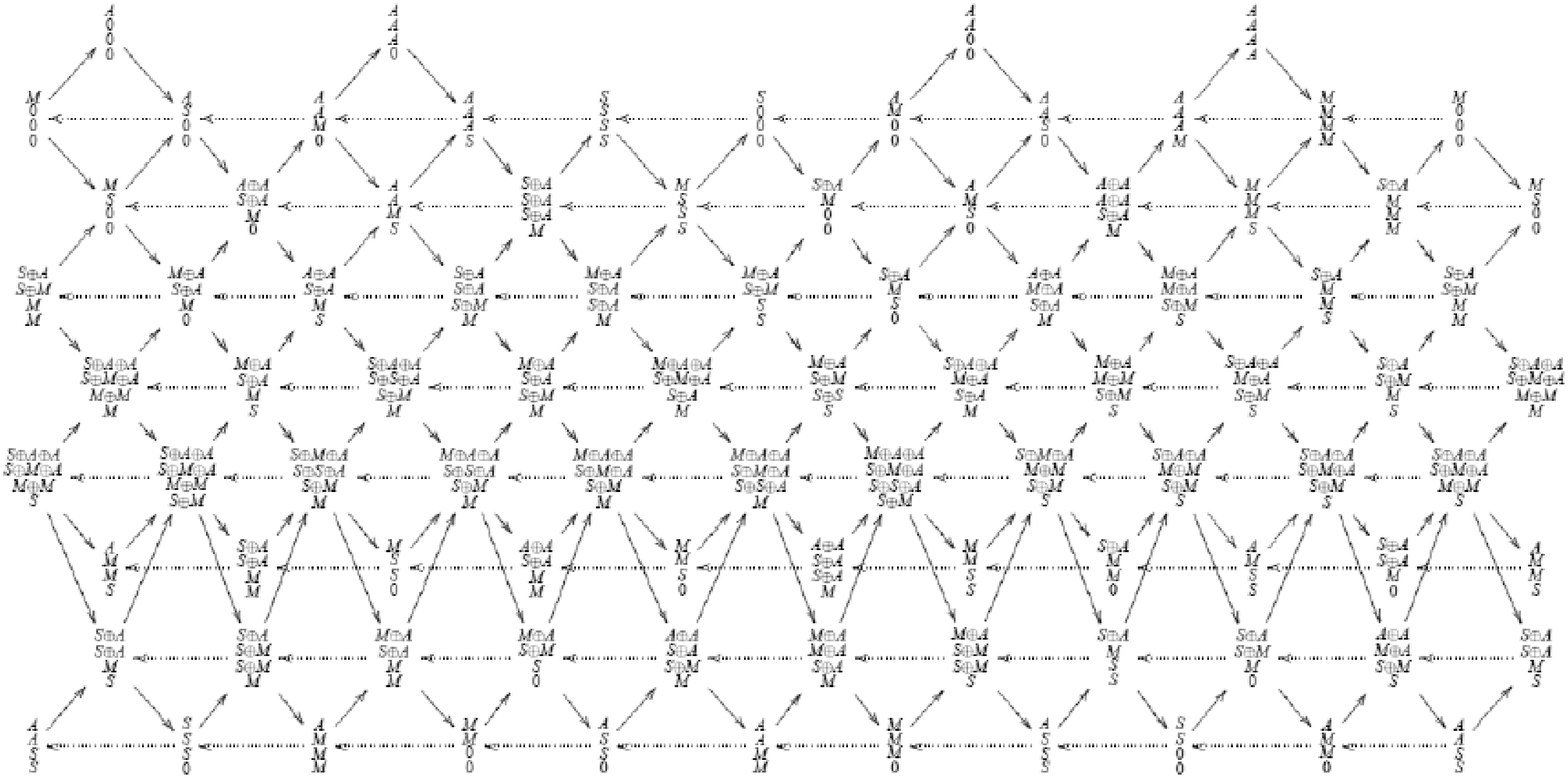}
\end{center}

\vskip10pt

\subsection{} Consider the selfinjective Nakayama algebra $\A = \A(2,
2)$, whose indecomposable modules are denoted by
$S_1=\xymatrix@=0.5cm{k\ar@<1ex>[r]&0\ar@<1ex>[l]}, \
S_2=\xymatrix@=0.5cm{0\ar@<1ex>[r]&k\ar@<1ex>[l]}, \
P_1=\xymatrix@=0.5cm{k\ar@<1ex>[r]^{1}&k\ar@<1ex>[l]^{0}}, \
P_2=\xymatrix@=0.5cm{k\ar@<1ex>[r]^{0}&k\ar@<1ex>[l]^{1}.}$ There
are 4 indecomposable projective objects and 6 indecomposable
nonprojective objects in $\mathcal{S}_2(\A)$, with the
Auslander-Reiten quiver as follows.

$$\xymatrix@R=0.1cm@C=0.3cm{& {\begin{smallmatrix}P_2\\0\end{smallmatrix}}\ar[dr] & &
{\begin{smallmatrix}P_2\\P_2\end{smallmatrix}}\ar[dr] & & & & \\
{\begin{smallmatrix}S_1\\0\end{smallmatrix}}\ar[ur]\ar[dr] & &
{\begin{smallmatrix}P_2\\S_1\end{smallmatrix}}\ar[ur]\ar[dr]\ar@{.>}[ll]
& &
{\begin{smallmatrix}S_2\\S_2\end{smallmatrix}}\ar[dr]\ar@{.>}[ll] &
&
{\begin{smallmatrix}S_1\\0\end{smallmatrix}}\ar[dr]\ar@{.>}[ll] & \\
& {\begin{smallmatrix}S_1\\S_1\end{smallmatrix}}\ar[ur] & &
{\begin{smallmatrix}S_2\\0\end{smallmatrix}}\ar[ur]\ar[dr]\ar@{.>}[ll]
& &
{\begin{smallmatrix}P_1\\S_2\end{smallmatrix}}\ar[ur]\ar[dr]\ar@{.>}[ll]
& &
{\begin{smallmatrix}S_1\\S_1\end{smallmatrix}}\ar@{.>}[ll]\\
& & & & {\begin{smallmatrix}P_1\\0\end{smallmatrix}}\ar[ur] & &
{\begin{smallmatrix}P_1\\P_1\end{smallmatrix}}\ar[ur] &}$$

There are 6 indecomposable projective objects and 12 indecomposable
nonprojective objects in $\mathcal{S}_3(\A)$, with the
Auslander-Reiten quiver as follows.
$$\xymatrix@R=0.08cm@C=0.4cm{& {\begin{smallmatrix}P_2\\0\\0\end{smallmatrix}}\ar[dr] & &
{\begin{smallmatrix}P_2\\P_2\\0\end{smallmatrix}}\ar[dr] & &
{\begin{smallmatrix}P_2\\P_2\\P_2\end{smallmatrix}}\ar[dr] & & & & &\\
{\begin{smallmatrix}S_1\\0\\0\end{smallmatrix}}\ar[dr]\ar[ur] & &
{\begin{smallmatrix}P_2\\S_1\\0\end{smallmatrix}}\ar[dr]\ar[ur]\ar@{.>}[ll]
& &
{\begin{smallmatrix}P_2\\P_2\\S_1\end{smallmatrix}}\ar[dr]\ar[ur]\ar@{.>}[ll]
& &
{\begin{smallmatrix}S_2\\S_2\\S_2\end{smallmatrix}}\ar[dr]\ar@{.>}[ll]
& &
{\begin{smallmatrix}S_1\\0\\0\end{smallmatrix}}\ar[dr]\ar@{.>}[ll] & &\\
& {\begin{smallmatrix}S_1\\S_1\\0\end{smallmatrix}}\ar[dr]\ar[ur] &
&
{\begin{smallmatrix}P_2\\S_1\\S_1\end{smallmatrix}}\ar[dr]\ar[ur]\ar@{.>}[ll]
& &
{\begin{smallmatrix}S_2\\S_2\\0\end{smallmatrix}}\ar[dr]\ar[ur]\ar@{.>}[ll]
& &
{\begin{smallmatrix}P_1\\S_2\\S_2\end{smallmatrix}}\ar[dr]\ar[ur]\ar@{.>}[ll]
& &
{\begin{smallmatrix}S_1\\S_1\\0\end{smallmatrix}}\ar[dr]\ar@{.>}[ll] & \\
& & {\begin{smallmatrix}S_1\\S_1\\S_1\end{smallmatrix}}\ar[ur] & &
{\begin{smallmatrix}S_2\\0\\0\end{smallmatrix}}\ar[dr]\ar[ur]\ar@{.>}[ll]
& &
{\begin{smallmatrix}P_1\\S_2\\0\end{smallmatrix}}\ar[dr]\ar[ur]\ar@{.>}[ll]
& &
{\begin{smallmatrix}P_1\\P_1\\S_2\end{smallmatrix}}\ar[dr]\ar[ur]\ar@{.>}[ll]
& &
{\begin{smallmatrix}S_1\\S_1\\S_1\end{smallmatrix}}\ar@{.>}[ll]\\
& & & & & {\begin{smallmatrix}P_1\\0\\0\end{smallmatrix}}\ar[ur] & &
{\begin{smallmatrix}P_1\\P_1\\0\end{smallmatrix}}\ar[ur] & &
{\begin{smallmatrix}P_1\\P_1\\P_1\end{smallmatrix}}\ar[ur] &}$$

\end{document}